\documentclass[smallextended,english]{svjour3}
\pdfoutput=1

\date{}
\def\makeheadbox{\relax}

\usepackage{xcolor}
\definecolor{lightgray}{gray}{0.75}

\usepackage[font=small]{caption}

\usepackage[OT1]{fontenc}
\usepackage[macce]{inputenc}
\usepackage{babel}
\usepackage{graphicx}
\usepackage{epstopdf}
\usepackage{amsmath}
\usepackage{amsfonts}
\usepackage{amssymb}
\usepackage{multirow}

\usepackage{tikz}
\newtheorem{thm}{Theorem}

\newtheorem{dfn}[thm]{Definition}

\begin{document}

\title{Parallelisation, initialisation, and boundary treatments for the diamond scheme}
\author{Stephen R Marsland
\and Robert I McLachlan
\and Matthew C Wilkins}
\institute{S R Marsland \at School of Mathematical and Computing Sciences, Victoria University of Wellington, New Zealand \and R  I McLachlan and M C Wilkins \at Institute of Fundamental Sciences, Massey University,
   Palmerston North, New Zealand}

\maketitle

\begin{abstract}\noindent
We study a class of general purpose linear multisymplectic integrators for Hamiltonian wave equations based on a diamond-shaped mesh. On each diamond, the PDE is discretized by a symplectic Runge--Kutta method. The scheme advances in time by filling in each diamond locally. We demonstrate that this leads to greater efficiency and parallelization and easier treatment of boundary conditions compared to methods based on rectangular meshes. We develop a variety of initial and boundary value treatments and present numerical evidence of their performance. In all cases, the observed order of convergence is equal to or greater than the number of stages of the underlying Runge--Kutta method.
\end{abstract}


\titlerunning{Parallelisation, initialisation and boundary treatments for the diamond scheme}

\section{Introduction}
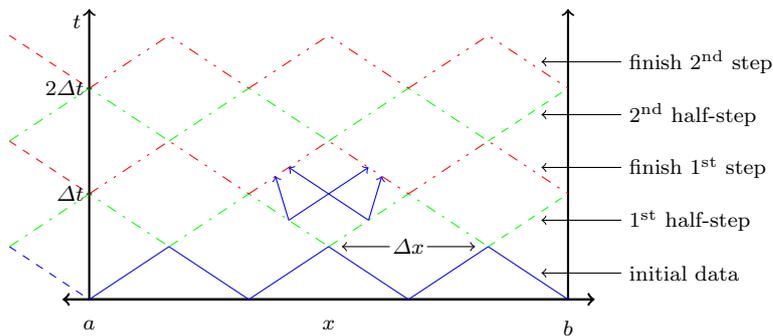
\begin{figure}%[htbp]
      \centering
      \begin{tikzpicture}[scale=0.35]

         % the domain
         \draw[<->,thick] (0, 11) -- (0, 0) -- (18, 0) -- (18, 11);
         \draw[<->,thick] (-1, 0) -- (19, 0);
         \draw[<->,thick] (-1, 0) -- (19, 0);
         
         % initial condition
         \draw[color=blue,style=dashed] (-3, 2) -- (0, 0);
         \foreach \x in {0cm, 6cm, 12cm} {
            \draw[color=blue] [xshift=\x] (0, 0) -- (3, 2) -- (6, 0);
         }
        
         % two steps
         \foreach \y in {0cm, 4cm} {
            \begin{scope}[yshift=\y]
               \foreach \x in {0cm, 6cm, 12cm} {
                  \draw[color=green,dash pattern=on 4pt off 3pt on 1.0pt off 3pt] [xshift=\x] (-3, 2) -- (0, 4) -- (3, 2);
               }
               \draw[color=green,style=dashed] (15, 2) -- (18, 4);
               \draw[color=red,style=dashed] (-3, 6) -- (0, 4);
               \foreach \x in {0cm, 6cm, 12cm} {
                  \draw[color=red,dash pattern=on 4pt off 3pt on 1.0pt off 3pt on 1.0pt off 3pt] [xshift=\x,yshift=4cm] (0, 0) -- (3, 2) -- (6, 0);
               }
            \end{scope}
         }
         
         % labels
         \draw (-3pt,4) -- (3pt,4) node[left] {$\Delta t$};
         \draw (0, 10.5) node[left] {$t$};
         \draw (-3pt,8) -- (3pt,8) node[left] {$2 \Delta t$};
         \draw (12,2) node {$\Delta x$};
         \draw[<-] (9.5,2) -- (11.3,2);
         \draw[->] (12.7,2) -- (14.5,2);
         \draw (0, -0.5) node[below] {$a$};
         \draw (9, -0.5) node[below] {$x$};
         \draw (18, -0.5) node[below] {$b$};
        
         % information
         \draw[->,blue] (7.5,3) -- (7.0,4.666);
         \draw[->,blue] (7.5,3) -- (10.5,5);
         \draw[->,blue] (10.5,3) -- (7.5,5);
         \draw[->,blue] (10.5,3) -- (11.0,4.666);

         % describe the levels
         \draw[<-] (17,1) -- (20,1) node[right] {initial data};
         \draw[<-] (17,3) -- (20,3) node[right] {$1^\mathrm{st}$ half-step};
         \draw[<-] (17,5) -- (20,5) node[right] {finish $1^\mathrm{st}$ step};
         \draw[<-] (17,7) -- (20,7) node[right] {$2^\mathrm{nd}$ half-step};
         \draw[<-] (17,9) -- (20,9) node[right] {finish $2^\mathrm{nd}$ step};
      
      \end{tikzpicture}

      \caption{Schematic of the diamond scheme for periodic boundary conditions.
      Information flows upwards as indicated by the solid blue 
      arrows for a typical diamond.  The solution, $\mathbf{z}$, is initialized on
      the solid blue zig-zag line.
      A step of
      the diamond scheme consists of two half-steps.  The first half-step
      calculates $\mathbf{z}$ along the green dash-dot line, which by periodicity is
      extended to the dashed line to the right.  The second half step uses
      the green dash-dot line to calculate the red dash-double-dotted line, which again by
      periodicity is extended to the left-hand dashed segment.}
      \label{fig:domain}
   \end{figure}

The diamond scheme is a family of fully discrete numerical methods for first-order hyperbolic PDEs
introduced in \cite{mclachlan2015multisymplectic}. It is based on the diamond grid shown in Figure \ref{fig:domain}.
The family is parameterised by its number of stages, $r$. The dependent variables
are associated with each of $r$ nodes on each edge of the grid; from data on the lower
two edges of a diamond, data on the top two edges can be computed {\em locally 
within a single diamond}. This feature is unique amongst schemes of such broad
applicability and motivates its further exploration. In particular, we argued in \cite{mclachlan2015multisymplectic}
that the local nature of the scheme indicated that it is highly parallelizable and
amenable to local boundary and initialization treatments. In this paper we present
numerical evidence supporting this view.

In particular, in Section 2 we present a serial and a parallel implementation of the diamond scheme for a nonlinear wave equation. The results are exceptionally good, showing high convergence orders and almost perfect
speedup with only $\gtrsim 5$ diamonds per processor. In Section 3 we develop two novel initialization methods and compare their performance to a reference method in which the diamonds are initialised with the exact solution of the differential equation. The observed orders of convergence are at least as good (and sometimes better) than those of the reference method. Dealing with
this unconventional aspect of the scheme, in which initial values are not determined by the problem data, is crucial. In Section 4 we develop local boundary treatments for nonhomogeneous Dirichlet and Neumann boundary conditions. These are tested on a variety of linear and nonlinear wave equations for extremely long integration times. In all the tests in all sections, the observed order
of convergence of the $r$-stage method is at least $r$, although often it exceeds this. Section 5 concludes.
The method is marked by its particular theoretical advantages and by its observed performance on a range of tests.

We describe the method as it applies to the family of Hamiltonian PDEs
\begin{equation}
   K\mathbf{z}_t + L\mathbf{z}_x = \nabla S(\mathbf{z}), \label{eq:hampde}
\end{equation}
where $K$ and $L$ are constant $n \times n$ real skew-symmetric matrices,
$\mathbf{z}\colon\Omega\to \mathbb{R}^n$, $\Omega\subset\mathbb{R}^2$, and $S\colon \mathbb{R}^n
\to\mathbb{R}$.  
Any  solutions to~\eqref{eq:hampde} satisfy the \emph{multisymplectic conservation law}
\begin{equation}
   \omega_t + \kappa_x = 0,    \label{eqn:symlaw}
\end{equation}
where $\omega = \frac{1}{2} (d\mathbf{z} \wedge K d\mathbf{z})$ and
$\kappa = \frac{1}{2} (d\mathbf{z} \wedge L d\mathbf{z})$~\cite{bridges2006numerical}\footnote{
That is, if ${\bf u}_1$, ${\bf u}_2$ are solutions to the variational equation
$K {\bf u}_t + L{\bf u}_x = S''({\bf z}){\bf u}$, then $\omega({\bf u}_1,{\bf u}_2) = {\bf u}_1^T K {\bf u}_2$
and $\kappa({\bf u}_1,{\bf u}_2) = {\bf u}_1^T L {\bf u}_2$.}.
A numerical method that satisfies a discrete version of
Eq.~\eqref{eqn:symlaw} is called a \emph{multisymplectic integrator};
see \cite{bridges2006numerical,simhamdyn} for reviews of multisymplectic integration.

For ODEs, there are effective symplectic integrators---such as symplectic Runge--Kutta methods---that apply to the entire class of ODEs $K z_t = \nabla S(z)$, and have excellent numerical properties, including
symplecticity, arbitrary order, small error constants, unconditional stability, and linear equivariance. 
One generalization of these methods to the PDEs \eqref{eq:hampde} is to apply high order Runge--Kutta methods in space and in time on a rectangular mesh \cite{reich}. This approach inherits some of the good features listed; some of its properties, including dispersion and order behaviour, are studied in \cite{mclachlan2014on}. However, the scheme does have some drawbacks. It is fully implicit, and it leads to discrete equations without a solution for periodic boundary conditions unless $r$ and $N$ (the number of cells in space) are {\em both} odd \cite{ryland,mclachlan2014on}. Solvability of the discrete equations is also affected by the boundary conditions, and no general effective treatment of boundary conditions is known for this method.

The first two issues, implicitness and boundary treatment, are related. They can be avoided for some PDEs, like the nonlinear wave equation, by applying suitably partitioned Runge--Kutta methods \cite{ryland,mclachlan2011linear,ryland2008multisymplecticity,ryland2007multisymplecticity}. When they apply, they lead to explicit ODEs amenable to explicit time-stepping, can have high order, and can deal with general boundary conditions. However, the partitioning means that they are not linear methods.

The diamond scheme is a different generalization of symplectic Runge--Kutta methods from ODEs to PDEs. It provides an approach that is multisymplectic, applies to all PDEs of the form (\ref{eq:hampde}), is linear in $z$, and is locally well-defined for any number of stages. 
We first give the definition of the class of diamond schemes that we will consider.

\begin{dfn}\label{def:diamond}
 A {\em diamond scheme} for the PDE~\eqref{eq:hampde} is a quadrilateral mesh in space-time together with a mapping of each quadrilateral  to a square to which a Runge--Kutta method is applied in each dimension, together with
initial data specified at sufficient edge points such that the solution can be propagated forward in time by locally solving for pairs of adjacent edges. 
\end{dfn}

It is convenient to first map 
each diamond in Figure~\ref{fig:domain} to a unit square using the linear transformation $T$ defined by (omitting unnecessary additive constants)
\begin{equation} \label{eqn:xttilde}
T\colon\quad   \tilde{x} = \tfrac{1}{\Delta x} x + \tfrac{1}{\Delta t} t, \quad \tilde{t} = -\tfrac{1}{\Delta x} x + \tfrac{1}{\Delta t} t.
\end{equation}
By the chain rule, Eq.~\eqref{eq:hampde} transforms to
\begin{equation}
\label{eq:hampdemod}
   \tilde{K} \tilde{\mathbf{z}}_{\tilde{t}} + \tilde{L} \tilde{\mathbf{z}}_{\tilde{x}} = \nabla S(\tilde{\mathbf{z}}),
\end{equation}
where
\begin{equation} \label{eqn:KLtilde}
   \tilde{K} = \tfrac{1}{\Delta t} K - \tfrac{1}{\Delta x} L, \quad \tilde{L} = \tfrac{1}{\Delta t} K + \tfrac{1}{\Delta x} L,
\end{equation}
and $\tilde z(\tilde x,\tilde t)=z(x,t)$.
As the PDE \eqref{eq:hampdemod} is still of the same class as \eqref{eq:hampde}, we may apply
the multisymplectic Runge-Kutta
collocation method given by Reich~\cite{reich} to Eq. (\ref{eq:hampdemod}) within a single unit square;
applying $T^{-1}$ yields a method on the diamond lattice.

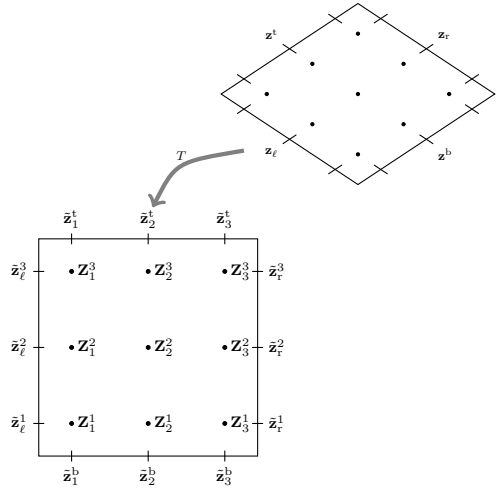
\begin{figure}%[htb]
      \centering
      \begin{tikzpicture}[scale=0.6,transform shape]

         % the diamond
         \draw (-3, 0) -- (0, -2) -- (3, 0) -- (0, 2) -- cycle;
%         \draw (0, 2) node[above]{$\mathbf{z}_0^1$};
%         \draw (0, -2) node[below]{$\mathbf{z}_0^{-1}$};
%         \draw (3, 0) node[right]{$\mathbf{z}_1^{0}$};
%         \draw (-3, 0) node[left]{$\mathbf{z}_{-1}^0$};

         % labels
         \draw [xshift=0.5cm,yshift=-1.667cm] (-0.15, 0.1) -- (0.15, -0.1); 
         \draw [xshift=1.5cm,yshift=-1cm] (-0.15, 0.1) -- (0.15, -0.1) node[below right]{$\mathbf{z}^\mathrm{b}$};
         \draw [xshift=2.5cm,yshift=-0.333cm] (-0.15, 0.1) -- (0.15, -0.1); 

         \draw [xshift=-2.5cm,yshift=-0.333cm] (-0.15, -0.1) -- (0.15, 0.1); 
         \draw [xshift=-1.5cm,yshift=-1cm]  (-0.15, -0.1) node[below left]{$\mathbf{z}_\ell$} -- (0.15, 0.1);
         \draw [xshift=-0.5cm,yshift=-1.667cm] (-0.15, -0.1) -- (0.15, 0.1); 

         \begin{scope}[xshift=-3cm,yshift=2cm]
         \draw [xshift=0.5cm,yshift=-1.667cm] (-0.15, 0.1) -- (0.15, -0.1); 
         \draw [xshift=1.5cm,yshift=-1cm] (-0.15, 0.1) node[above left]{$\mathbf{z}^\mathrm{t}$} -- (0.15, -0.1);
         \draw [xshift=2.5cm,yshift=-0.333cm] (-0.15, 0.1) -- (0.15, -0.1); 
         \end{scope}

         \begin{scope}[xshift=3cm,yshift=2cm]
         \draw [xshift=-2.5cm,yshift=-0.333cm] (-0.15, -0.1) -- (0.15, 0.1); 
         \draw [xshift=-1.5cm,yshift=-1cm]  (-0.15, -0.1) -- (0.15, 0.1)   node[above right]{$\mathbf{z}_\mathrm{r}$};
         \draw [xshift=-0.5cm,yshift=-1.667cm] (-0.15, -0.1) -- (0.15, 0.1); 
         \end{scope}

         % internal points
         \filldraw   (-2, 0) circle (1pt) 
         (-1, 0.667) circle (1pt) 
         (0, 1.333) circle (1pt) 
         (-1, -0.667) circle (1pt) 
         (0, 0) circle (1pt) 
         (1, 0.667) circle (1pt) 
         (0, -1.333) circle (1pt) 
         (1, -0.667) circle (1pt) 
         (2, 0) circle (1pt) ;

         % the arrow down to the square
         \draw[->,gray,ultra thick] (-2.5,-1.25) .. controls (-4.0,-1.5) .. (-4.5, -2.5) node[pos=0.5,above,black] {$T$};

         \begin{scope}[xshift=-7cm, yshift=-8cm, scale=1.2]

         % the square
         \draw (0, 0) -- (4, 0) -- (4, 4) -- (0, 4) -- cycle;
         %\draw[<->] (0, 4.5) -- (0, 0) -- (4.5, 0);

         % labels
%         \foreach \x / \xtext in {0, 0.6/1, 2, 3.4/3}
         \foreach \x / \xtext in {0.6/1, 2, 3.4/3}
            \draw [xshift=\x cm,yshift=0cm] (0, 0.1) -- (0, -0.1) node[below]{$\tilde{\mathbf{z}}^\mathrm{b}_\xtext$};
         \foreach \x / \xtext in {00.6/1, 2, 3.4/3}
            \draw [xshift=\x cm,yshift=4cm] (0, 0.1) node[above]{$\tilde{\mathbf{z}}^\mathrm{t}_\xtext$} -- (0, -0.1);
         \foreach \y / \ytext in {00.6/1, 2, 3.4/3}
            \draw [xshift=0cm,yshift=\y cm] (-0.1, 0)  node[left]{$\tilde{\mathbf{z}}_\ell^\ytext$} -- (0.1, 0);
         \foreach \y / \ytext in {0.6/1, 2, 3.4/3}
            \draw [xshift=4cm,yshift=\y cm] (-0.1, 0) -- (0.1, 0) node[right]{$\tilde{\mathbf{z}}_\mathrm{r}^\ytext$};
         \foreach \y / \ytext in {0.6/1, 2, 3.4/3} {
            \foreach \x / \xtext in {0.6/1, 2, 3.4/3} {
               \filldraw (\x, \y) circle (1pt) node[right]{$\mathbf{Z}_\xtext^\ytext$};
            }
         }
            \end{scope}

      \end{tikzpicture}

      \caption{The diamond transformed by a linear transformation,
      $T$, to the unit square.  The square
      contains $r \times r$ ($r=3$ in this example) internal stages, $\mathbf{Z}_i^j$.
      The solution is known along the bottom and left hand sides.
      The method proceeds as two sets of $r$ Gauss Runge--Kutta $r$-stage
      methods: internal stage values, $\mathbf{Z}_i^j,\ \mathbf{X}_i^j,\ \mathbf{T}_i^j$, are calculated, then the
      right and top updated.}
      \label{fig:single_diamond}
   \end{figure}

Let $(A,b)$ be the parameters of an $r$-stage Runge--Kutta method. In what
follows, we will take the method to be the Gauss Runge--Kutta method.
Figure~\ref{fig:single_diamond} shows a
diamond with $r = 3$, and its transformation to the unit square.
The square contains $r \times
r$ internal grid points, as determined by the Runge-Kutta coefficients $c$,
 and internal stages $\mathbf{Z}_i^j$, which are analogous to the usual
internal grid points and stages in a Runge-Kutta method. The internal stages also carry 
the variables $\mathbf{X}_i^j$ and $\mathbf{T}_i^j$ which approximate $z_x$ and $z_t$, respectively, at the internal stages.

The dependent variables of the method are the values of $z$ at the edge grid points.
%Let $I=\left\{1, \ldots, r\right\}$ so that, for example, $\tilde{\mathbf{z}}^\mathrm{b}_I$ refers
%to $\tilde{\mathbf{z}}^\mathrm{b}_i, i = 1, \ldots, r $.
%If the  qualifier $I$ does not appear
%then the left- or bottom-most
%corner is included, for example $\tilde{\mathbf{z}}^\mathrm{b}$ refers to 
%$\tilde{\mathbf{z}}^\mathrm{b}_i, i = 0, \ldots, r $.
%Note $\tilde{\mathbf{z}}^\mathrm{b}_0 = \tilde{\mathbf{z}}_\ell^0=\mathbf{z}_0^{-1}$.
The Runge--Kutta discretization is
\begin{align}
   \mathbf{Z}_i^j &= \tilde{\mathbf{z}}_\ell^j + \sum_{k=1}^r a_{ik} \mathbf{X}_k^j, \label{eqn:msZ1} \\
   \mathbf{Z}_i^j &= \tilde{\mathbf{z}}^\mathrm{b}_i + \sum_{k=1}^r a_{jk} \mathbf{T}_i^k, \label{eqn:msZ2} \\
   \nabla S(\mathbf{Z}_i^j) &= \tilde{K} \mathbf{T}_i^j + \tilde{L} \mathbf{X}_i^j, \label{eqn:msZ3}
\end{align}
together with the update equations
\begin{align}
   \tilde{\mathbf{z}}_\mathrm{r}^i &= \tilde{\mathbf{z}}_\ell^i + \sum_{k=1}^r b_k  \mathbf{X}_k^i, \label{eqn:msupdate1} \\
   \tilde{\mathbf{z}}^\mathrm{t}_i &= \tilde{\mathbf{z}}^\mathrm{b}_i + \sum_{k=1}^r b_k\mathbf{T}_i^k, \label{eqn:msupdate2}
\end{align}
for $i,j =1,\dots,r$.  Here
$\tilde{\mathbf{z}}_\ell^i$ and
$\tilde{\mathbf{z}}^\mathrm{b}_i$ are known.
Eqs.~\eqref{eqn:msZ1}--\eqref{eqn:msZ3} are first solved for the
internal stage values $\mathbf{Z}_i^j$, $\mathbf{X}_i^j$, and $\mathbf{T}_i^j$, then Eqs.~\eqref{eqn:msupdate1}
and~\eqref{eqn:msupdate2} are used to calculate 
$\tilde{\mathbf{z}}^\mathrm{t}_i$ and
$\tilde{\mathbf{z}}_\mathrm{r}^i$.  
Eqs.~\eqref{eqn:msZ1}--\eqref{eqn:msZ3} 
are $3r^2$
equations in $3r^2$ unknowns 
$\mathbf{Z}$, $\mathbf{X}$, and $\mathbf{T}$.
Eqs.~\eqref{eqn:msZ1} and~\eqref{eqn:msZ2} are linear in 
$\mathbf{X}$ and $\mathbf{T}$. 
Thus in practice the method requires solving a set of $r^2n$ nonlinear
equations for $\mathbf{Z}$ in each diamond.

The method does not use values at the corners. If necessary,
solutions at the corners can be obtained using Runge--Kutta
update equations along the edges, combined with averages
where two edges meet.

The following basic properties of the diamond scheme
are established in \cite{mclachlan2015multisymplectic}. The conservation law in Theorem \ref{thm:diamond_conservation_law}
is a  discretization of the integral of $\omega_t+\kappa_x=0$ over a single
diamond, transferred to the boundary of the diamond using Stokes's theorem
and discretized by Gauss quadrature.

\begin{thm} \label{thm:diamond_conservation_law}
The diamond scheme satisfies the discrete symplectic conservation law
$$
\frac{1}{\Delta t} \sum_{i=1}^r b_i ( \omega_i^\mathrm{t} + \omega_\mathrm{r}^i - (\omega_\ell^i + \omega_i^\mathrm{b}) ) +
\frac{1}{\Delta x} \sum_{i=1}^r b_i ( \kappa_\mathrm{r}^i + \kappa_i^\mathrm{b} - (\kappa_i^\mathrm{t} + \kappa_\ell^i))  = 0,
$$
where 
$\omega_n^m = \tfrac{1}{2} d \mathbf{z}_n^m \wedge K d \mathbf{z}_n^m$,
$\kappa_n^m = \tfrac{1}{2} d \mathbf{z}_n^m \wedge L d \mathbf{z}_n^m$.
\end{thm}
%(Refer to Figure~\eqref{fig:single_diamond}
%for the definition of these labels).

We also recall the result of \cite{mclachlan2014on} that the Reich collocation scheme 
with the $r$-stage Gauss method,
when applied on a rectangular mesh to the hyperbolic PDE \eqref{eq:hampde}
with initial conditions on $t=0$ and periodic boundary conditions, 
has global errors of order at least $r$. In some cases,
order $r+1$ is observed, which is the stage order of the method. 
Therefore we expect convergence of order at least $r$
from the diamond scheme in cases where it is stable.

One major difference between the diamond and rectangular meshes is that
the diamond mesh means that the method is effectively multi-stage. Indeed,
it is the extra initial data, at different time levels, that allows the diamonds
to be filled in independently. In the ODE case $L=0$, the diamond scheme
yields an $r$-step integrator whose underlying 1-step method is the original
Gauss method. In this paper we test the diamond scheme numerically
on several different  wave equations $u_{tt} - u_{xx} = f(u)$
on the domain $x\in[a,b]$, $t\in[0,T]$. The initial conditions are
Cauchy (i.e., $u(x,0)$ and $u_t(x,0)$ are given), and various
boundary conditions (periodic, Dirichlet, and Neumann) are applied at $x=a$ and $x=b$.
Note that there are other methods (in particular,
Lobatto IIIA--IIIB in space and an explicit symplectic splitting method in time)
that perform outstandingly well on nonlinear wave equations. We adopt this
equation simply as a first test: if the method fails here it is almost certainly fails overall.

The wave equation has several formulations of the form \eqref{eq:hampde}. Here 
we use the formulation with $v=u_t$, $w=u_x$, $\mathbf{z} = (u, v, w)^T$, 
\begin{gather} \label{eqn:1dwaveeqn}
   K = \begin{pmatrix} 0 & -1 & 0\\1 & 0 & 0\\0 & 0 & 0 \end{pmatrix},
   \quad L = \begin{pmatrix} 0 & 0 & 1\\0 & 0 & 0\\-1 & 0 & 0
   \end{pmatrix}, \\ \quad S(\mathbf{z}) = -V(u)+\frac{1}{2}v^2 - \frac{1}{2}w^2, \;\mathrm{and}\;
   f(u)=V'(u).
\end{gather} 
Then we have the following stability result.

\begin{thm} \cite{mclachlan2014on} 
   The  diamond scheme with $r=1$ applied to the wave equation $u_{tt}=u_{xx}$ with periodic
   boundary conditions  is linearly stable
   when
   $\lambda = \frac{\Delta t}{\Delta x} \le 1$.
\end{thm}

%We present convergence and parallelisation results in Section \ref{chapter:diamondimplem},
%describe two successful methods for imposing initial conditions in Section \ref{sec:diamondinit},
%and methods for dealing with Dirichlet and Neumann boundary conditions
%in Section \ref{sec:bcs}. Numerical experiments on several different initial--boundary
%value problems and forcings indicate that the methods are stable and accurate.

\section{Diamond implementation} \label{chapter:diamondimplem}

Three implementations of the diamond scheme were prepared:
serial implementations in Python and C and a parallel version
in C. Although the scheme is most naturally expressed using
rank-3 tensors, for the implementations the tensors
were flattened to matrices. Solving the nonlinear
equations on each diamond requires a nonlinear solver.  
The Python code used the
SciPy~\cite{scipy} routine \verb+fsolve+, which is a wrapper around the
\textsc{minpack}~\cite{minpack} \verb+hydrd+ and \verb+hybrj+
algorithms, which are based on the Powell hybrid
method~\cite{Powell01011964}.  The C codes used the
GNU Scientific Library~\cite{gsl} routines \verb+gsl_multiroot_fsolver_hybrids+ and
\verb+gsl_multiroot_fdfsolver_gnewton+ which again are wrappers around
the \textsc{minpack}  \verb+hydrd+ and \verb+hybrj+ algorithms.
The \verb+hydrd+ algorithm approximates the Jacobian, whereas
\verb+hybrj+ requires the exact Jacobian; both versions were
implemented and compared.

\subsection{Convergence of the serial implementation}
\begin{figure}
   \centering
   \includegraphics[width=9cm]{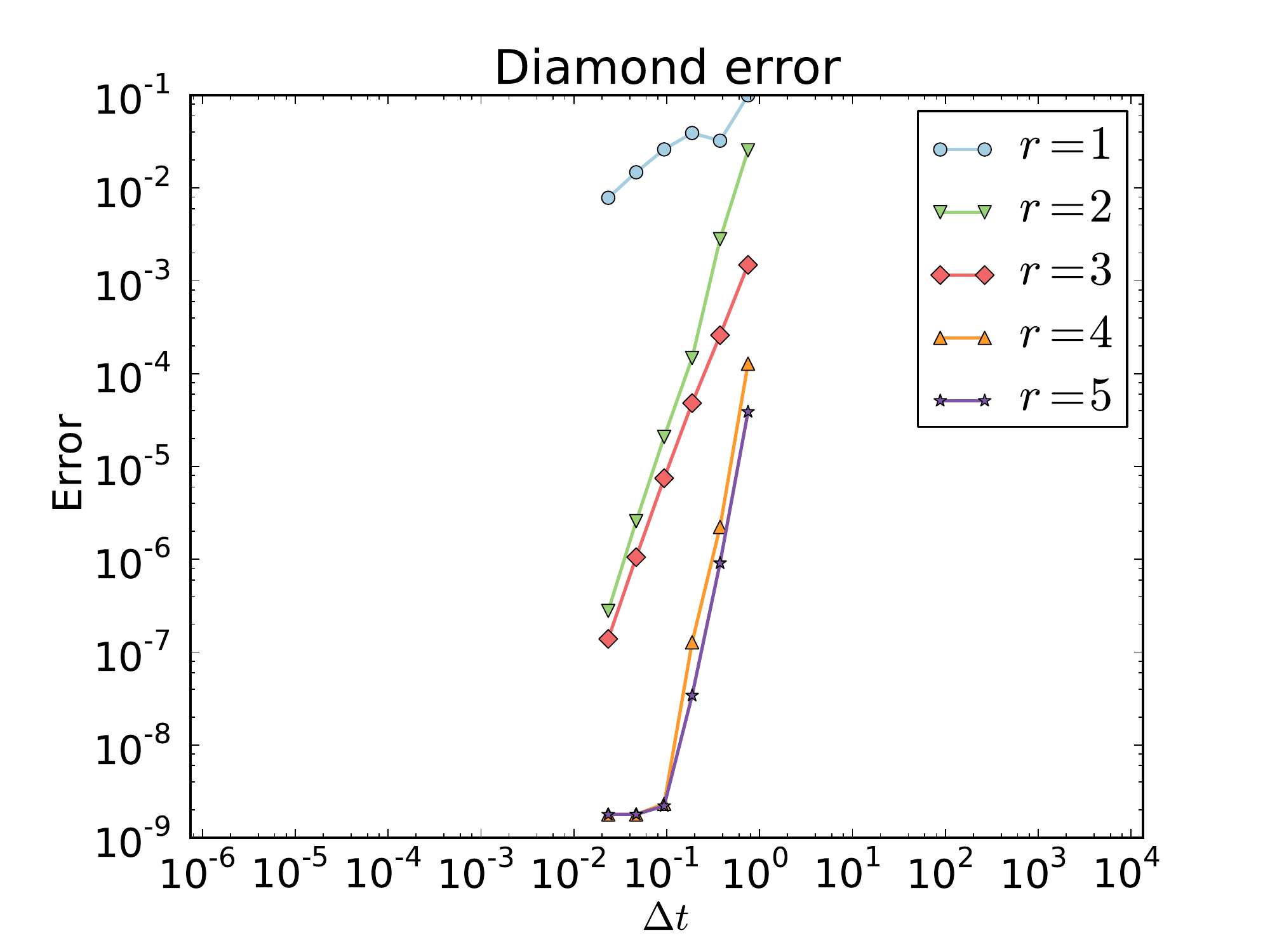}
\caption{The error of the diamond scheme with varying $r$ applied to
the
multi-symplectic Hamiltonian PDE arising from the sine-Gordon
equation.
The true solution was the so-called \emph{breather} on the domain
$\left[-30, 30\right]$.
The Courant number is fixed at
$\frac{1}{2}$ as $\Delta t$ is decreased.  The order of the method
appears to be $r$ when $r$ is odd and $r+1$ when $r$ is even.}

   \label{fig:error_r12345_sinegordon_periodic_s2_internal0}
\end{figure}

We first test the convergence of the serial implementation on 
a nonlinear wave equation. 
The diamond scheme with varying $r$ was used to solve the sine-Gordon
equation  $u_{tt}-u_{xx}=-\sin(u)$.
An exact solution is the so-called \emph{breather},
\begin{equation} \label{eqn:breather}
u(x,t) = 4 \arctan \left( \frac{\sin\left(\frac{t}{\sqrt{2}}\right)}{ \cosh\left(\frac{x}{\sqrt{2}}\right)} \right).
\end{equation}
The domain is taken significantly large, $[a,b]=\left[-30, 30\right]$, so
the solution can be assumed to be periodic. 
The error is the discrete 2-norm of $u$,
\begin{equation}
\label{eq:E}
E^2 = \frac{b - a}{N} \sum_i^N \left(u_i - u(a + i \Delta x, T) \right)^2.
\end{equation}
The number of diamonds at each time level is $N=40,80,\dots,1280$, and the integration time, $T$, is twice the largest time step. 
The Courant number $\frac{\Delta t}{\Delta x} = \frac{1}{2}$ is held fixed.
The $2rN$ initial values of $z=(u,u_t,u_x)$ needed at the bottom edge of the
first row of diamonds are provided by the exact solution.
The results for the global error are shown in
Fig.~\eqref{fig:error_r12345_sinegordon_periodic_s2_internal0}.  It is apparent that for this problem the
method converges, and the
order
appears to be $r$ when $r$ is odd and $r+1$ when $r$ is even.

\subsection{Parallel implementation and speedup}

In the parallel implementation of the diamond scheme,
the domain is divided into strips, finite in width in the $x$
direction and potentially infinite in the $t$ direction.  The width $L$
of the domain is  a multiple of $\Delta x$.  There are $N = \frac{L}{\Delta x}$ diamonds
in each row, and $p$ processors.  The $N$ diamonds are divided as equally as
possible into $p$ contiguous regions.  If $N$ is a multiple of $p$
then each processor will get an equal number of diamonds to work on.
Otherwise $p - n$ processors get $k$ diamonds, and $n$ processors get
$k+1$ diamonds, where $n = N-pk$ and $k = \lfloor \frac{N}{p} \rfloor$.

Each processor calculates the solution on its strip of diamonds.  After
the first half time-step (see \ref{fig:domain}) the solution at the right edge of the
strip must be passed to the right-hand neighbour and received
from the left.  This results in $p$ transmits of
vectors of length  $n(r+1)$ in length.  At the next half time-step the
left hand edge must be passed to the left and received from the right,
another $p$ transmits.  If the solution is to be output each
processor must send a number of values---twice the number of diamonds
in its sub-domain---to one processor.  This communication was achieved
with an MPI gather.  Since this processor has very slightly more work to perform it is
desirable to ensure it is one of the $p-n$ processors that receives
$k$ diamonds to work on.  
 
To test how much the diamond scheme could benefit from running in
parallel it was run on the New Zealand eScience Infrastructure's
(NeSI's) Pan
Intel Linux cluster physically located at the University of Auckland,
New Zealand.  At the time of use  the cluster had approximately
6000 cores each running somewhere between 2GHz and 3GHz with most
around 2.7GHz or 2.8GHz.  Most cores had at least 10GB of RAM
available, far more than we required.  Due to the
busyness of the cluster it was impractical to request specific CPUs
for each job, thus there was a certain variability in timing tests
simply because of the different speeds of the CPUs available.
However the uni-processor jobs, arguably the most important while
testing parallel speed-ups, did run on the most common 2.7GHz or
2.8GHz processors.

The diamond scheme was
initialized with the diamond initialization method detailed in
Section~\ref{sec:diamondinit}, $r$ was set to 5, the number of
time steps was 1000, $\Delta t=0.05$, and the periodic sine-Gordon problem
from Table~\ref{table:sampleperiodicproblems} was used.  Despite the
cluster having approximately 6000 cores, by trial and error it was
apparent that only about a maximum of 300 or 400 cores could be readily
available on demand.  So each trial consisted of nine runs with
the number of cores being 1, 3, 7, 20, 56, 100, 150, 300, and 350.
For each run the wall-clock time was recorded using the Unix
\texttt{date} command, the program run, and then the wall-clock time
checked again.  Each trial (set of nine runs) was performed twice with a couple of days in
between each trial, and the time results averaged.

\begin{figure}
   \centering
   \includegraphics[width=5.6cm]{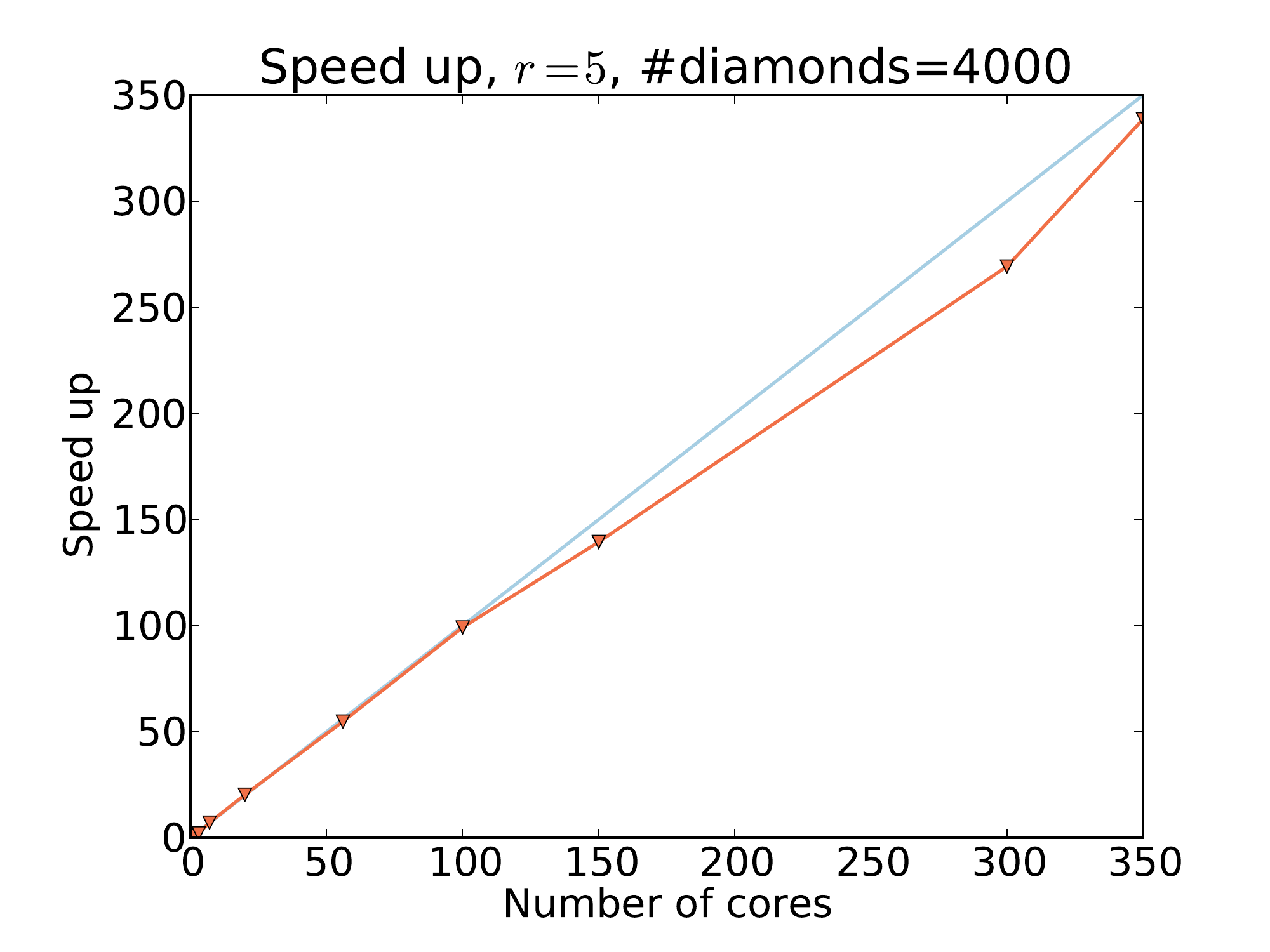}
   \includegraphics[width=5.6cm]{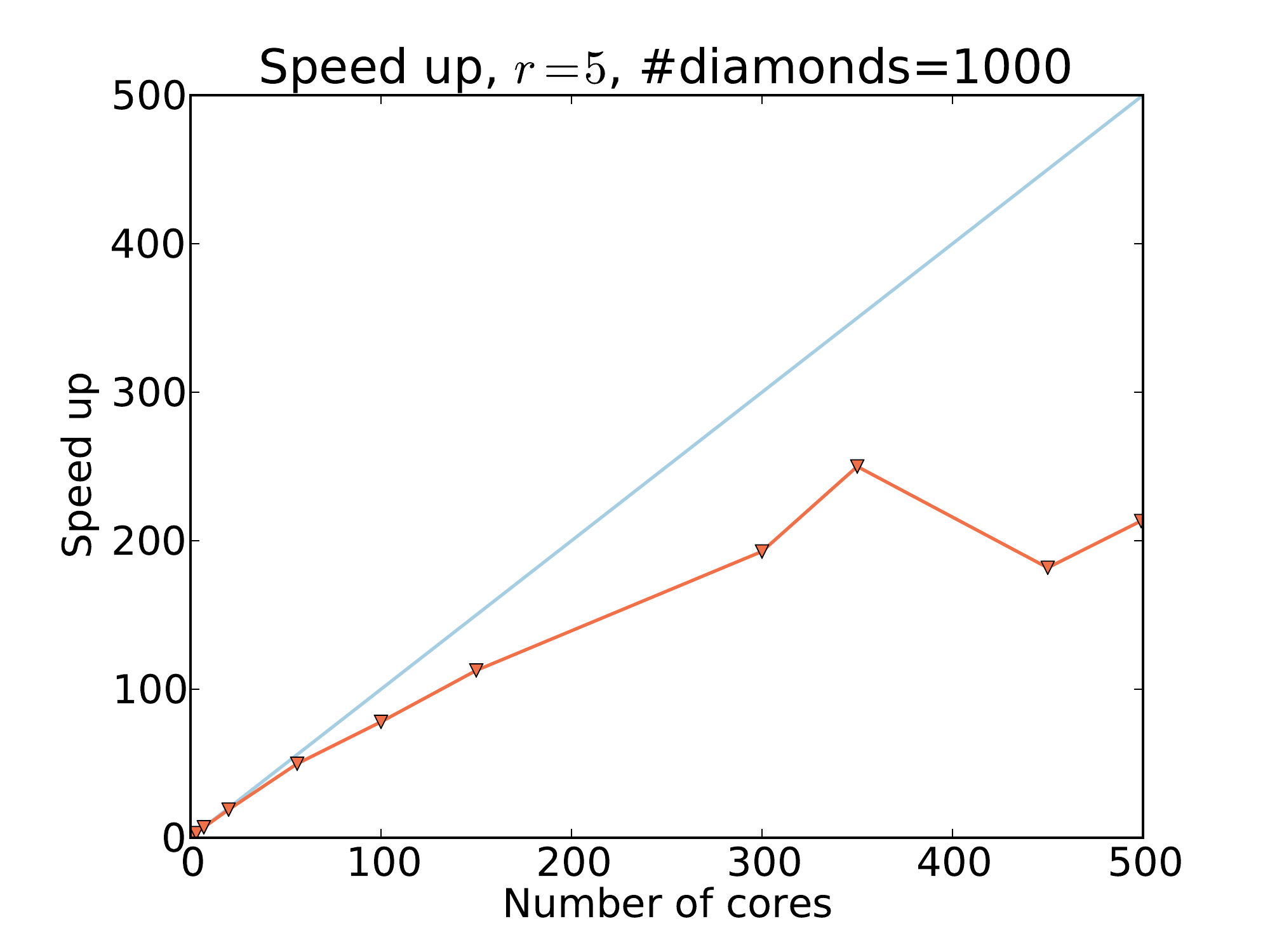}
   \caption{Speed-up of the diamond scheme versus the number of
   cores for the code running on the NeSI Pan cluster.  Code that was perfectly
   parallelizeable would have the speed-up equal to the number of
cores (the blue line).  As the ratio of the number of cores to the
amount of work (number of diamonds across the domain) increases one
would expect the speed-up to deviate from the perfect blue line. 
On the left, with at least $4000/350>11$ diamonds per core,
the speedup is very good. On the right, with as few as 2 diamonds 
per core, the performance deteriorates. The speed up was calculated from a single run (not an
average of two runs) for the 450 and 499 number of cores runs.
\label{fig:mpi:init=diamond:physics=per_sinegordon:r=5:s=2}
}
\end{figure}

According to Amdahl's law~\cite{amdahl}, for a particular problem
size, if $n$ is the number of cores, and $B \in
[0,1]$ is the fraction of the algorithm that is strictly serial,
then the theoretical time $T(n)$ for the algorithm to run on $n$ cores
is
$$
   T(n) = T(1)\left(B + \frac{1}{n}(1-B)\right).
$$
Thus the theoretical speed-up $S(n)$ is
\begin{equation} \label{eqn:speedup}
S(n) = \frac{T(1)}{T(n)} = \frac{1}{B + \tfrac{1}{n}(1-B)}.
\end{equation}
Letting $n \rightarrow \infty$ gives a theoretical maximum speed-up of
$\tfrac{1}{B}$.   By increasing $n$ until the speed-up begins to
tail-off it is possible to estimate $B$.  
For a perfectly parallelizeable algorithm the speed-up should be equal
to the number of cores used.  In the first trial, shown on the left in
Figure~\ref{fig:mpi:init=diamond:physics=per_sinegordon:r=5:s=2},
$\Delta x$ was such that there were 4000 diamonds across the domain.

As the ratio of the number of cores to the amount of work (number of
diamonds across the domain) increases, one would expect the speed-up
to diverge from the perfect speed-up line.  This is because the overhead in
communication will gradually swamp the gains in computation time.  For this
trial the speed-up is still very good and it is impossible to
estimate $B$, the fraction of the algorithm that is strictly serial.
Ideally, the number of cores
would be increased until a deviation from the perfectly parallelizable
line could be reliably detected, however no more cores were easily
available.

So, instead of increasing the number of cores, the number of diamonds
was decreased.
Figure~\ref{fig:mpi:init=diamond:physics=per_sinegordon:r=5:s=2} (right)
shows the second trial where the number of diamonds across the domain
was decreased to 1000.  One of the trials included two extra runs with
$n=450$ and $n=499$.  Because it took many days for these runs to
begin executing, the second trial did not include these large runs,
and no averaging could take place.  This figure shows the speed-up
reaching approximately 250 before beginning to tail off.  So for this
size of problem, from~\eqref{eqn:speedup} this equates to $B \approx
0.4\%$, which is remarkably low.  
The conclusion is that the diamond scheme is exceptionally
parallelizable.

\section{Diamond scheme initialization} \label{sec:diamondinit}

The diamond mesh creates an issue for initialization which is not
present (or hardly present) for rectangular meshes.
In this section we develop two initialization methods,
the {\em diamond initialization} and {\em phantom initialization} methods, and
compare them to a reference `method' in which
initial values are taken from the exact solution.

\subsection{Diamond initialization method}

For the nonlinear wave equations, initial conditions for $u$ and $u_t$ are specified
at $t=0$. Differentiating with respect to $x$, we can assume that
$u_x$, and hence all components of $z$, are known at $t=0$. 

The $x$-axis cuts the first row of diamonds
in half, yielding a row
of triangles.  In the diamond initialization method, each triangle is mapped to
the unit square, then the usual set of
equations~\eqref{eqn:msZ1}--~\eqref{eqn:msupdate2} solved, giving
values for $z$ on the top-left and top-right edges of the triangle.

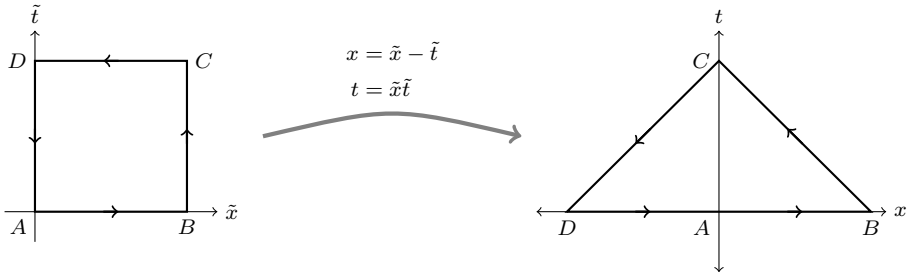
\begin{figure}
   \centering
   \begin{tikzpicture}[scale=2]

      \draw[->] (-0.2, 0) -- (1.2, 0) node[right] {$\tilde{x}$};
      \draw[->] (0, -0.2) -- (0, 1.2) node[above] {$\tilde{t}$};
      \draw[thick] (0,0) -- (1, 0) -- (1, 1) -- (0, 1) -- cycle;
      \draw[->,thick] (0.45, 0) -- (0.55,0);
      \draw[->,thick] (1, 0.45) -- (1, 0.55);
      \draw[->,thick] (0.55, 1) -- (0.45,1);
      \draw[->,thick] (0, 0.55) -- (0, 0.45);
      \draw (0,0) node[below left] {$A$};
      \draw (1,0) node[below] {$B$};
      \draw (1,1) node[right] {$C$};
      \draw (0,1) node[left] {$D$};

      \draw[->,gray,ultra thick] (1.5,0.5) .. controls (2.35,0.7) .. (3.2, 0.5) node[pos=0.5,above,black] {
\begin{minipage}{2cm}\begin{align*}
x &= \tilde{x} - \tilde{t}\\
t &= \tilde{x}\tilde{t}
\end{align*}\end{minipage}};

      \begin{scope}[xshift=4.5cm,yshift=0.0cm]
      \draw[<->] (-1.2, 0) -- (1.1, 0) node[right] {$x$};
      \draw[<->] (0, -0.4) -- (0, 1.2) node[above] {$t$};
      \draw[thick] (-1,0) node[below] {$D$} -- (0, 0) node[below left] {$A$} -- (1, 0) node[below] {$B$} -- (0, 1) node[left] {$C$} -- cycle;
      \draw[->,thick] (-0.55, 0) -- (-0.45,0);
      \draw[->,thick] (0.45, 0) -- (0.55,0);
      \draw[->,thick] (0.55, 0.45) -- (0.45, 0.55);
      \draw[->,thick] (-0.45, 0.55) -- (-0.55,0.45);

      \end{scope}

   \end{tikzpicture}

   \caption{The unit square under the map $x = \tilde{x} - \tilde{t}, t = \tilde{x}\tilde{t}$.}
   \label{fig:transform_unitsquare_to_triangle}
\end{figure}

The transformation, $x = \tilde{x} - \tilde{t}, t =
\tilde{x}\tilde{t}$,
illustrated in Figure~\ref{fig:transform_unitsquare_to_triangle},
takes the unit square to the triangle $(-1,0),(1,0),(0,1)$.

Adding a translation and scaling results in the map
\begin{align} \label{eqn:unitsquare_to_triangle}
   x &= \frac{\Delta x}{2}(\tilde{x} - \tilde{t}) + b \\
   t &= \frac{\Delta t}{2}\tilde{x} \tilde{t}
\end{align}
which takes the unit square to the triangle $(b-\tfrac{\Delta x}{2},0), (b+\tfrac{\Delta x}{2},0),(b,\tfrac{\Delta t}{2})$.
Recall that the transformed $K$ and $L$ are given by
\begin{align*}
   \tilde{K} &= g_t  K + g_x  L \\
   \tilde{L} &= f_t  K + f_x  L,
\end{align*}
where $(\tilde{x}, \tilde{t}) = (f(x,t),g(x,t))$.  
This yields
\begin{align*}
   \tilde{K} &= \tfrac{2}{\Delta x \Delta t (\tilde{x}+\tilde{t})} \left( \Delta x  K - \tilde{t} \Delta t  L \right) \\
   \tilde{L} &= \tfrac{2}{\Delta x \Delta t (\tilde{x}+\tilde{t})} \left( \Delta x  K + \tilde{x} \Delta t  L \right).
\end{align*}

For initializing the diamond scheme, values of $z$ are needed on the
bottom zig-zag (Figure~\ref{fig:domain}) spaced
according to the Runge--Kutta vector $c$.  Because the above
map~\eqref{eqn:unitsquare_to_triangle} and its
inverse are linear on the edges, this same spacing can be used in
$(\tilde{x}, \tilde{t})$ space.

\begin{table}
   \centering
   \begin{tabular}{c c c}
      \hline
      $r$ & \multicolumn{2}{c}{Order} \\
        & exact & diamond \\ \hline 
      1 & 1 & 2\\
      2 & 3 & 3\\
      3 & 3 & 5\\
      4 & 5 & 5\\
      5 & 5 & 6\\ \hline
   \end{tabular}
   \caption{Observed order of convergence 
      of the diamond scheme initialized with the
   exact solution and with the diamond scheme initialization. (Data
   obtained from Fig.~\ref{fig:error_r12345_periodic_sinegordon_s2_diamond_initexact_init}.)  It is
apparent that for this problem the diamond initialization performs as
well, or better, than the exact initialization.  The order appears to
be $r+1$ for most $r$ (for $r=4$ the order is $r+2$ for some reason),
whereas for the exact initialization the order is $r$ ($r$ odd) and
$r+1$ ($r$ even).}
   \label{table:orderdiamondinit}
\end{table}

\begin{figure}
   \centering
   \includegraphics[width=9cm]{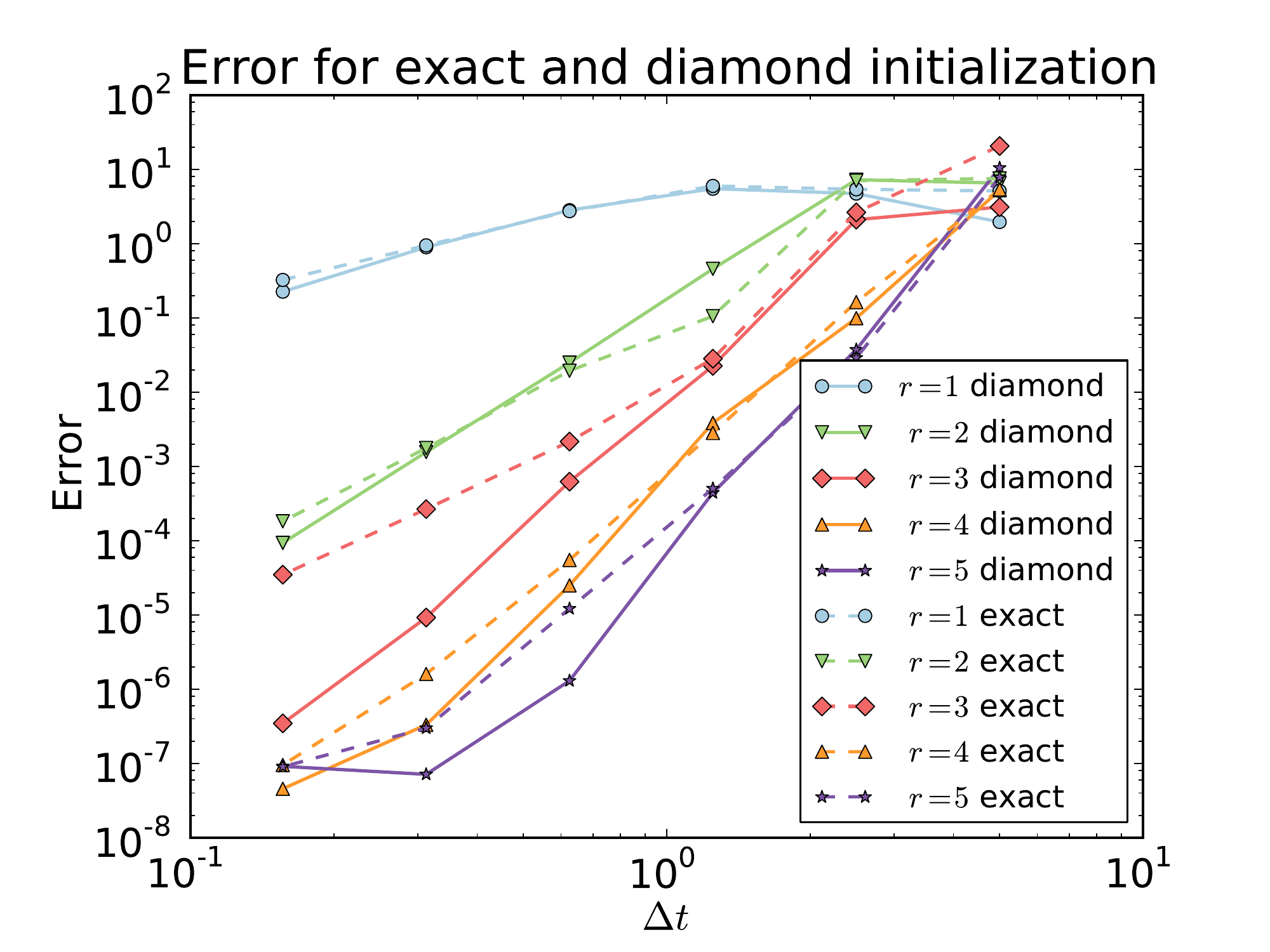}
   \caption{The error of the diamond scheme initialized using the
      exact and diamond methods applied to
the sine-Gordon equation.
The exact solution is a so-called \emph{breather} on  $\left[-30,
30\right]$.  The Courant number is fixed at $\frac{1}{2}$ as $\Delta t$ is
decreased.  The error is the 2-norm of the global error at the final time, Eq. \ref{eq:E}.
For this problem the diamond initialization is as accurate or
better than exact initialization.
}
   \label{fig:error_r12345_periodic_sinegordon_s2_diamond_initexact_init}
\end{figure}

Figure~\ref{fig:error_r12345_periodic_sinegordon_s2_diamond_initexact_init}  shows the error
of exact and diamond initialization as $\Delta t$ is reduced while keeping the Courant
number $\frac{\Delta t}{\Delta x} = \frac{1}{2}$.  The integration time $T=1.5$
is twice the largest time step.  It is apparent that for this problem
the diamond initialization is equal, or better, than exact
initialization.  Table~\ref{table:orderdiamondinit} shows the
observed convergence order of the two initialization methods.

\begin{table}
   \centering
   \begin{tabular}{l l l l}
      \hline
      Name & Equation & Range & Solution \\ \hline
      Esin     & $u_{tt} + u_{xx} = 0$ & $0 \le x \le 2 \pi$ & $e^{2 \sin(x-t-3)}$ \\
      Sincos   & $u_{tt} + u_{xx} = 0$ & $0 \le x \le 2 \pi$ & $\sin(x)\cos(t)$ \\
      Coscos   & $u_{tt} + u_{xx} = -u$ & $0 \le x \le \pi$ & $\cos(2x)\cos(\sqrt{5}t)$ \\
      sine-Gordon   & $u_{tt} + u_{xx} = -\sin(u)$ & $-30 \le x \le 30$ & $4 \arctan \left( \frac{\sin\left(\frac{t}{\sqrt{2}}\right)}{ \cosh\left(\frac{x}{\sqrt{2}}\right)} \right)$ \\ \hline
   \end{tabular}
   \caption{Sample problems.}
   \label{table:sampleperiodicproblems}
\end{table}
We next apply the diamond initialization method to the four different
wave equations with different forcings and initial conditions
given in Table~\ref{table:sampleperiodicproblems}.
In each case the boundary conditions are periodic.
The number of diamonds at each time level is $N=10,20,\dots,1280$, and the integration time, $T=1.5$ is twice the largest time step.

The computed global  errors are shown in
Figure~\ref{fig:error:init=diamond:physics=per_expasin:r=12345:s=2}.
Table~\ref{table:erroroffourproblems} shows the observed convergence order of the
diamond scheme for these problems.  It is apparent that for this problem, the
order is at least $r$.  

\def\fw{5.6cm}
\begin{figure}%[htbp]
   \centering
   \includegraphics[width=\fw]{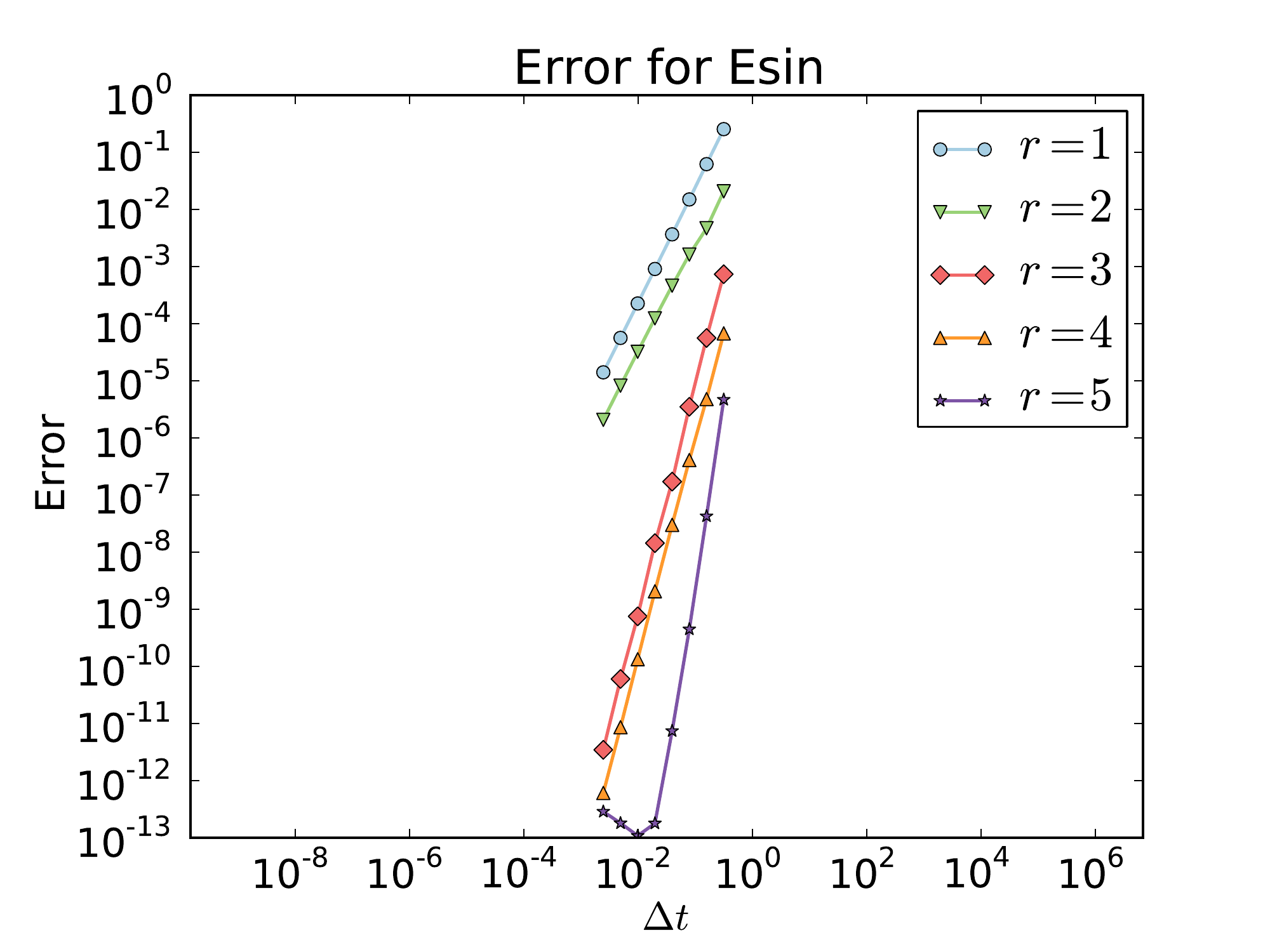}
   \includegraphics[width=\fw]{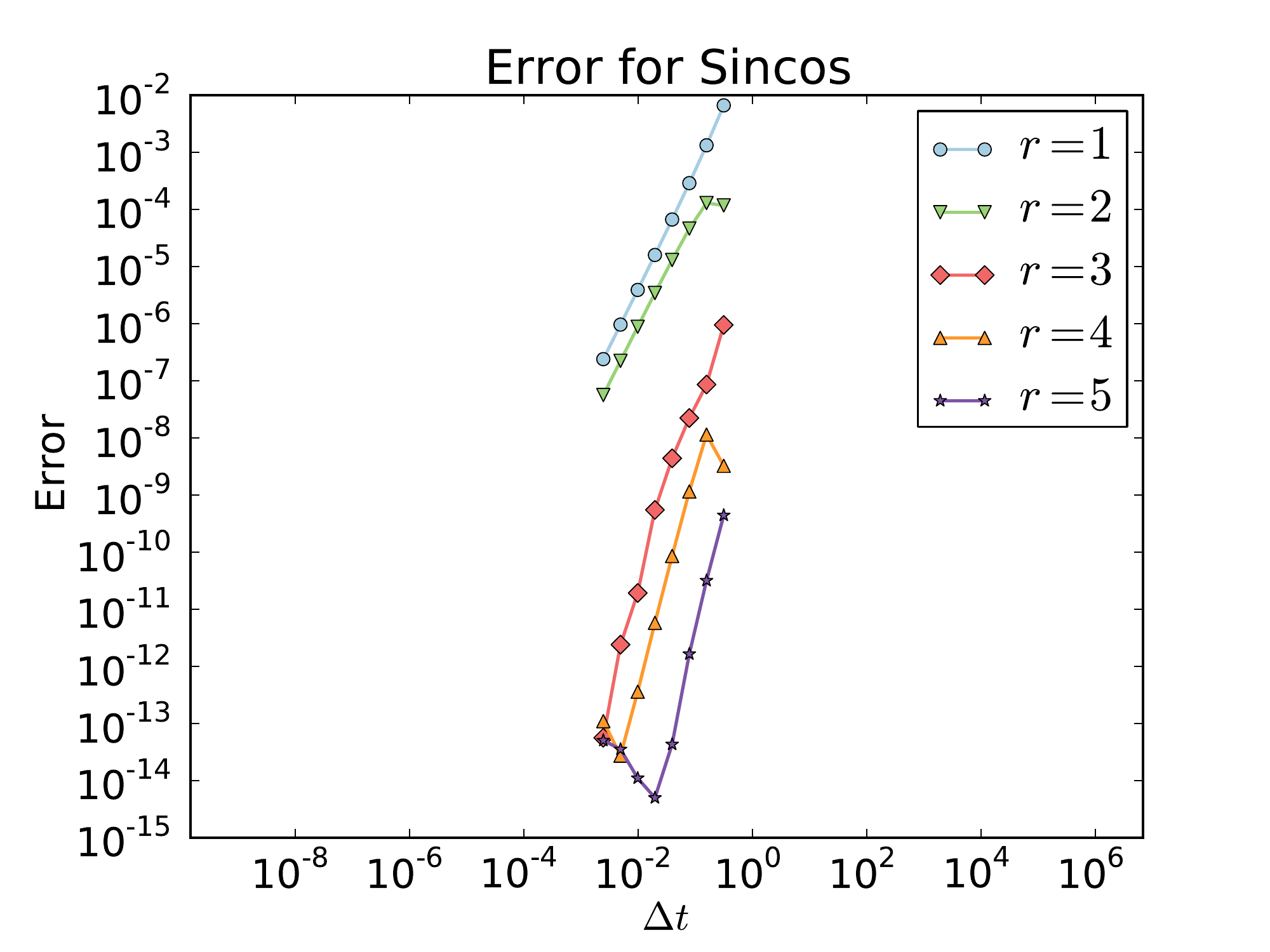}
   \includegraphics[width=\fw]{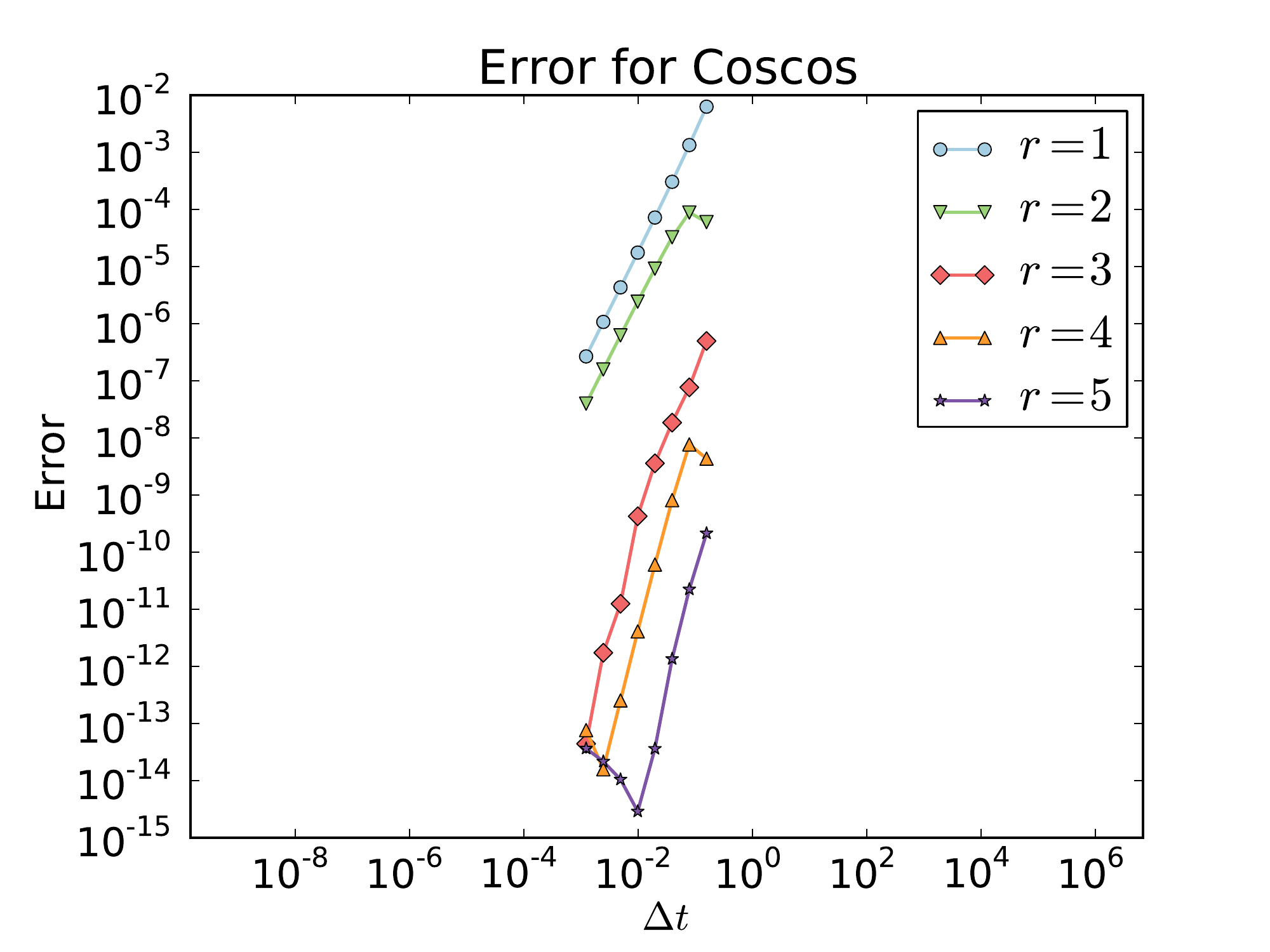}
   \includegraphics[width=\fw]{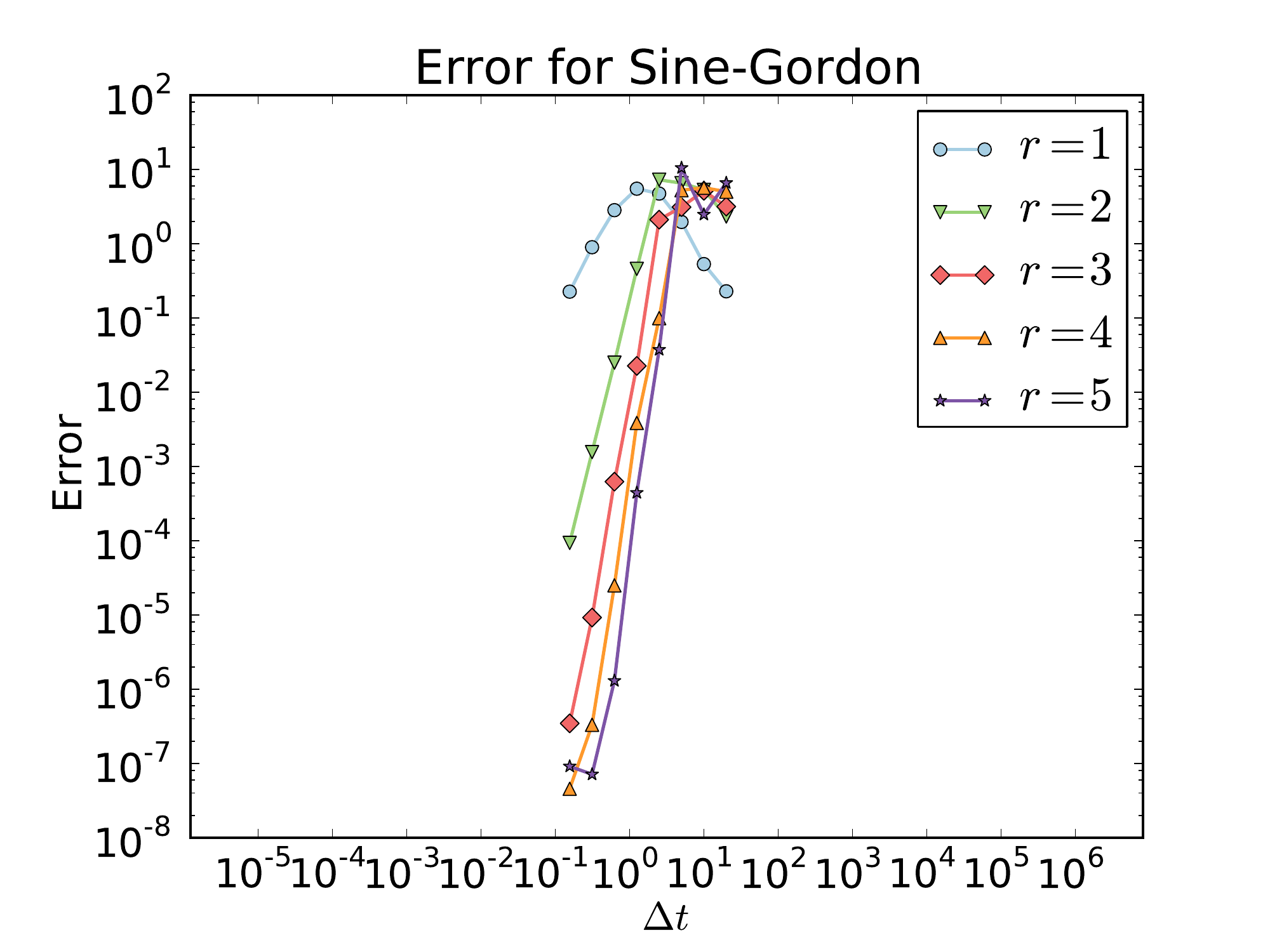}
   \caption{The error of the diamond scheme with varying number of stages $r$ applied
      to the Esin, Sincos, Coscos, and sine-Gordon problems (see
      Table~\ref{table:sampleperiodicproblems}).  The Courant number
      is fixed at $\frac{1}{2}$ as $\Delta t$ is decreased.
      Table~\ref{table:erroroffourproblems} summarizes the observed convergence
      order given by
       the slope of these lines.     \label{fig:error:init=diamond:physics=per_expasin:r=12345:s=2}
      }
\end{figure}

\begin{table}
   \centering
   \begin{tabular}{c c c c c}
      \hline
      $r$ & \multicolumn{4}{c}{Order} \\
          & Esin & Sincos & Coscos & sine-Gordon \\ \hline
      1 & 2.0 & 2.1 & 2.1 & 1.8 \\
      2 &1.9 & 1.9 & 1.9 & 4.1 \\
      3 & 4.0 & 3.4 & 3.3 & 5.6 \\
      4 & 3.8 & 3.7 & 3.8 & 5.4 \\
      5 & 6.2 & 4.8 & 4.6 & 6.8 \\ \hline
   \end{tabular}
   \caption{Observed convergence order
      of the problems given in
      Table~\ref{table:sampleperiodicproblems} (data from
      Figure~\ref{fig:error:init=diamond:physics=per_expasin:r=12345:s=2}). 
      The observed convergence order of the $r$ stage method is at least $r$.
      }
   \label{table:erroroffourproblems}
\end{table}

\subsection{Boundary initialization method}

\begin{figure}
      \centering
      \begin{tikzpicture}[scale=0.85]
         \draw[dashed] (-3, 0) -- (0, -2) -- (3, 0);
         \draw (-3, 0) -- (0, 2) -- (3, 0);
         \draw (-3.5, 0) -- (3.5, 0);
         \draw (0.75, -1.5) node {$\bigstar$};
         \draw (2.25, -0.5) node {$\bigstar$};
         \draw (-0.75, -1.5) node {$\bigstar$};
         \draw (-2.25, -0.5) node {$\bigstar$};

         \draw[dashed] (1.5, 0) circle (5pt);
         \draw[dashed] (-1.5, 0) circle (5pt);
         \draw (0, -1.0) node {$\bigstar$};
         \draw (0, 1.0) node {$\bigstar$};
      \end{tikzpicture}
      \caption{An initial phantom diamond in the $r=2$
         scheme.  The solution, $z$, is not known on the SW or SE
         edges.
         Internally,
         nothing is known at the stars, and some information is known
         at the dashed circles.  There needs to be enough known at the
         dashed circles to match the missing information
         on the SW and SE stars.}
      \label{fig:r2dirichletinit}
   \end{figure}
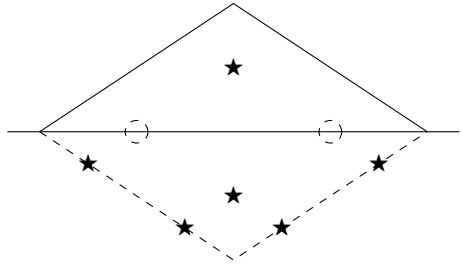

This method is inspired by the successful treatment
of Dirichlet and Neumann boundary conditions (described
in Section \ref{sec:bcs}) that utilises a phantom diamond.
 Here a phantom diamond is
constructed about the $t=0$ axis as illustrated in Figure~\ref{fig:r2dirichletinit}.

Compared to an internal diamond, the SW and SE edges are
missing values of $z$; these will now become  free variables,
whose values will be determined by the Runge--Kutta update equations.
To compensate, more information must be
gathered from the initial conditions at $t=0$. Note that $r$ of the $r^2$ internal
stages lie on $t=0$ and hence have partial data available.

By not specifying the values of $2nr$ values of $z$ on the SW and SE
edges of the phantom diamond, at additional $2nr$ values must
be specified at the internal stages on $t=0$. 
For the wave equation, $n=3$, so 6 values are
required per internal stage on $t=0$ axis.  The initial
conditions specify $u$ and $v$ ($=u_t$), giving
two components of $Z$ ($\approx z$) at $t=0$. 
The remaining required data can be obtained by differentiating
the PDE. 
Differentiating in $x$
gives $w = u_x$, $v_x = u_{tx}$, and $w_x=u_{xx}$. This specifies
the value of two components of $X$ ($\approx z_x$) 
at the $t=0$ internal stages. Because $u_{tx} =
u_{xt}$, $w_t$ is also known; this specifies
the value of one component of $T$ ($\approx z_t$).
In total we now have 6  out of the 9 components
of $Z$, $X$, and $T$ specified at each $t=0$ internal
stage. The Runge--Kutta equations now
yield a closed system that can be solved locally
within each phantom diamond separately, and
the update equations yield the values of $z$
on the NW and NE edges.

\begin{figure}
   \centering
   \includegraphics[width=9cm]{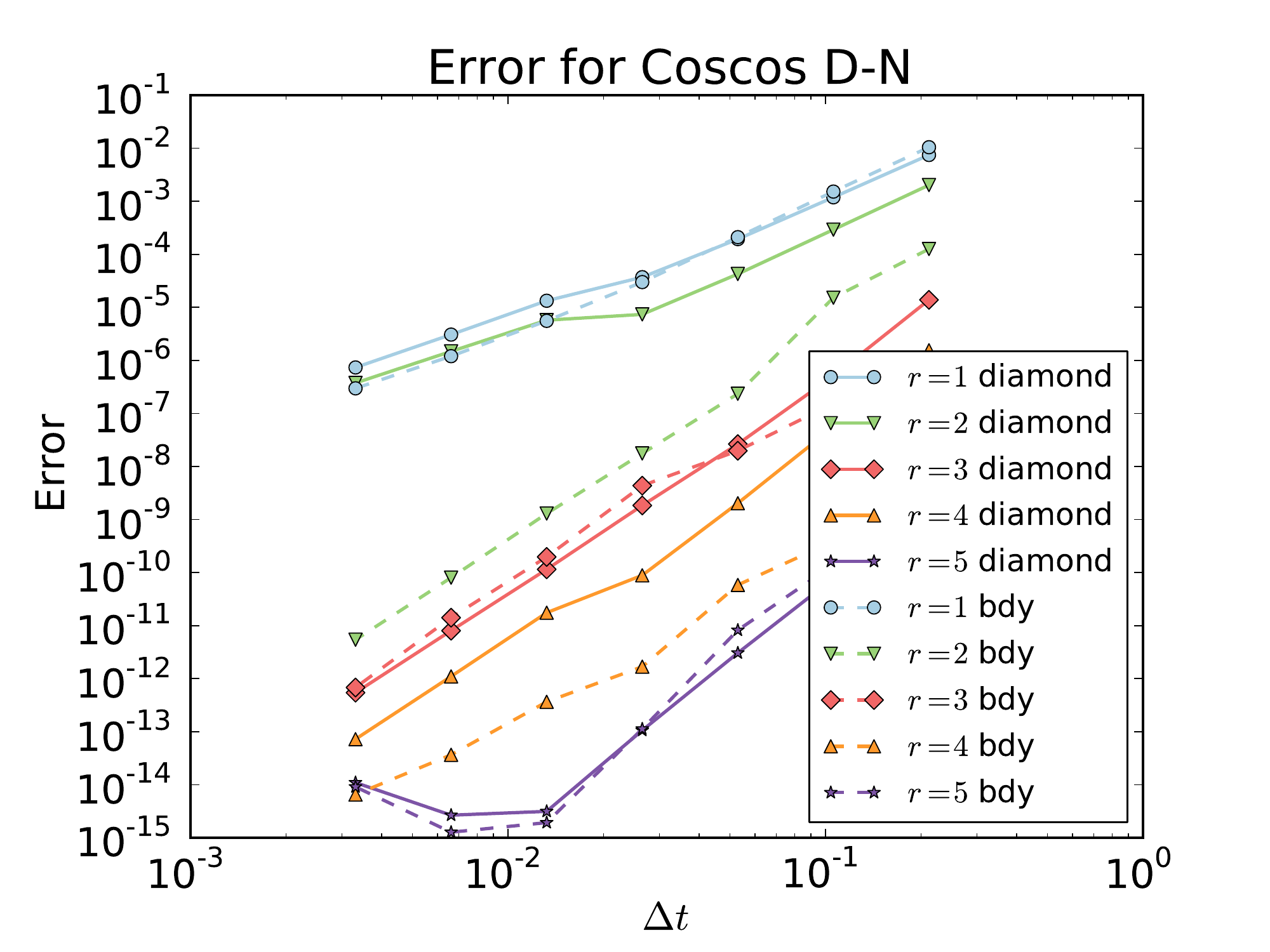}
   \caption{The error of the diamond scheme initialized with the diamond method
      (Section~\ref{sec:diamondinit}) and the boundary initialization
      method (Section~\ref{sec:bcs}), with varying $r$
      applied to the Coscos D-N problem (see
      Table~\ref{table:samplednproblems}).  The Courant number is
      fixed at $\frac{1}{2}$ as $\Delta t$ is decreased.  The boundary
      initialization method is as good as, or better, than the diamond
      initialization method.}
   \label{fig:error:init=diamondANDbdy:physics=d_n_cossquared:r=12345:s=2}
\end{figure}

The boundary method does not appear to adversely affect the
order of the scheme: it is at least $r$.
Figure~\ref{fig:error:init=diamondANDbdy:physics=d_n_cossquared:r=12345:s=2}
shows the error for the Coscos D-N problem initialized with the usual
diamond method, and the boundary method.  It is apparent that the
boundary initialization method does as well, or better, than the diamond
initialization.  

Although the boundary initialization method works well in this
example, it is a little ad-hoc as it relies on being able
to compute enough information from the initial conditions.
This may depend on the PDE and its formulation.

\section{Dirichlet and Neumann boundary conditions} \label{sec:bcs}

The construction of stable, high-order methods for hyperbolic
initial--boundary values problems is not easy. To cite one
successful approach,
a significant develop effort over many years has resulted in stable
finite difference methods using the \emph{summation by parts}
and \emph{simultaneous approximation term}
methods~\cite{kreiss1977existence,2014JCoPh.268...17S,carpenter1994time,olsson1995summation}.
These finite difference operators approximate $u_x$ (resp. $u_{xx}$) at all points, using
different finite differences near the boundary. Stability is achieved by requiring that the finite
difference is skew- (resp. self-) adjoint with respect to an inner product, designed
along with the method. In comparison, the compactness of the diamond scheme
indicates that we might hope to construct entirely local boundary treatments, systematically
for all $r$.

   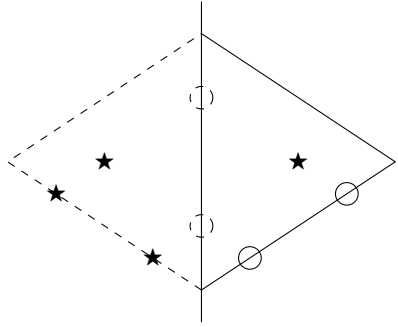
\begin{figure}
      \centering
      \begin{tikzpicture}[scale=0.85]
         \draw[dashed] (0, -2) -- (-3, 0) -- (0, 2);
         \draw (0, 2) -- (3, 0) -- (0, -2);
         \draw (0, -2.5) -- (0, 2.5);
         \draw (0.75, -1.5) circle (5pt);
         \draw (2.25, -0.5) circle (5pt);
         \draw (-0.75, -1.5) node {$\bigstar$};
         \draw (-2.25, -0.5) node {$\bigstar$};

         \draw (1.5, 0) node {$\bigstar$};
         \draw (-1.5, 0) node {$\bigstar$};
         \draw[dashed] (0, -1.0) circle (5pt);
         \draw[dashed] (0, 1.0) circle (5pt);
      \end{tikzpicture}
      \caption{A left-hand boundary phantom diamond in the $r=2$
         scheme.  The solution, $z$, is known at the circles on the SE
         edge.  For an internal diamond, $z$ would also be known at the
         stars on the SW edge, but not in this case.  Internally,
         nothing is known at the stars, and some information is known
         at the dashed circles.  There needs to be enough known at the
         dashed circles to match the missing information
         on the SW stars.}
      \label{fig:r2dirichlet}
   \end{figure}

At the left and right boundaries, the geometry of the grid alters; the
diamonds are cut in half to become triangles.
Suppose a boundary condition specifies $k$ values of $z$.
For the wave equation with Dirichlet boundary conditions,
$n=3$ and $k=1$ if $u(a,t)=g(t)$ is given.
Inspecting the whole diamond in Figure~\ref{fig:r2dirichlet},
$nr$ data values are missing on the SW edge.
These need to be made up from the boundary conditions.
Typically, just imposing the $k$ boundary conditions
at the $r$ internal stages on $x=a$ is not sufficient
to get a closed system.

We have developed and tested many approaches. A basic
requirement is that the resulting method should be solvable
and stable subject to a CFL condition, ideally $\Delta t / \Delta x < 1$. 
In addition, we found that some methods that worked well on 
a simple test problem (e.g. the linear wave equation with
homogeneous Dirichlet boundary conditions) did not
work on more complicated problems. Thus, extensive
testing was required. Before describing the 
successful method and its behaviour, we list
some methods that were not robustly able to solve wave
equations with a variety of boundary conditions:
(i) specifying some components of $z$ at
more points on the boundary than just the
$r$ stages; (ii) using extra information from an adjacent
interior diamond; (iii) mapping the boundary triangle
to a square, as in the diamond initialisation method;
and (iv) combinations of these.

The method that was ultimately successful is
the following {\em boundary scheme} that
we describe for Dirichlet and Neumann boundary 
conditions for the nonlinear wave equation.
In both cases we use the phantom diamond
as shown in Fig.~\ref{fig:r2dirichlet}, with
$r$ conditions imposed at the
internal stages on $x=a$, to compensate
for the $r$ missing values at the SW edge.
The entire set of equations \eqref{eqn:msZ3} is
then solved simultaneously, and the NE edge
values filled in using the update equations
(\ref{eqn:msupdate1},\ref{eqn:msupdate2}). The  $r$ new conditions
are equations in the $3nr$ internal
dependent variables $Z$ ($\approx z$),
$X$ ($\approx z_x$), and $T$ $(\approx z_t)$.
These equations come from the boundary
conditions, their derivatives with
respect to $t$, and the derivatives
of the PDE with respect to $x$ and $t$.
However, note that we cannot use
the PDE itself as it is already imposed
at the internal stages.
\begin{enumerate}
\item For the Dirichlet boundary condition
$u(a,t)=g(t)$, at the $x=a$ stages
we specify the value of $u$ (1st component of
$Z$); differentiate the boundary condition
to get $u_t(a,t)=g'(t)$, and we 
specify $v(a,t)=g'(t)$ (2nd component of $Z$);
differentiating again gives $u_{tt}(a,t)=g'(t)$,
which together with the PDE gives
$w_x(a,t)=g''(t) - f'(g(t))$ (3rd component
of $X$).
\item
For the Neumann boundary condition
$u_x(a, t) = h(t)$, at the $x=a$ stages we specify $w$
(3rd component of $Z$); 
differentiating with respect to time gives $w_t=h'(t)$ (3rd
component of $T$); and by equality of mixed partial derivatives $v_x =
w_t$, so we specify $v_x=h'(t)$ (2nd component of $X$).
\end{enumerate}

\begin{table}
   \centering
   \begin{tabular}{l l l l}
      \hline
      Name & Domain & Left boundary  & Right boundary \\ \hline
      Esin DD       & $0.2 \le x \le \tfrac{\pi}{3}$ & Dirichlet & Dirichlet \\
      Sincos DD     & $0.2 \le x \le \tfrac{\pi}{3}$ & Dirichlet & Dirichlet \\
      Sincos DN     & $0.2 \le x \le \tfrac{\pi}{3}$ & Dirichlet & Neumann \\
      Coscos DD     & $0.2 \le x \le \tfrac{\pi}{3}$ & Dirichlet & Dirichlet \\
      Coscos DN     & $0.2 \le x \le \tfrac{\pi}{3}$ & Dirichlet & Neumann \\
      sine-Gordon  DD & $-2 \le x \le 2$ & Dirichlet & Dirichlet \\ \hline
   \end{tabular}
   \caption{Sample non-periodic problems.  
See Table~\ref{table:sampleperiodicproblems} for the exact equations and
solutions.  The Dirichlet boundary conditions are found using the exact
solution, and the Neumann conditions by differentiating the exact
solution with respect to $x$.
   \label{table:samplednproblems}}
\end{table}

This boundary scheme is now applied to the four sample
problems given in Table~\ref{table:sampleperiodicproblems} with a mix
of Dirichlet and Neumann boundary conditions.
Table~\ref{table:samplednproblems} summarizes the problems.
In each case $10^6$--$10^7$ time steps were computed to ensure
that the equations were solvable at each step and the solutions
remained bounded.

\begin{figure}
   \centering
   \includegraphics[width=\fw]{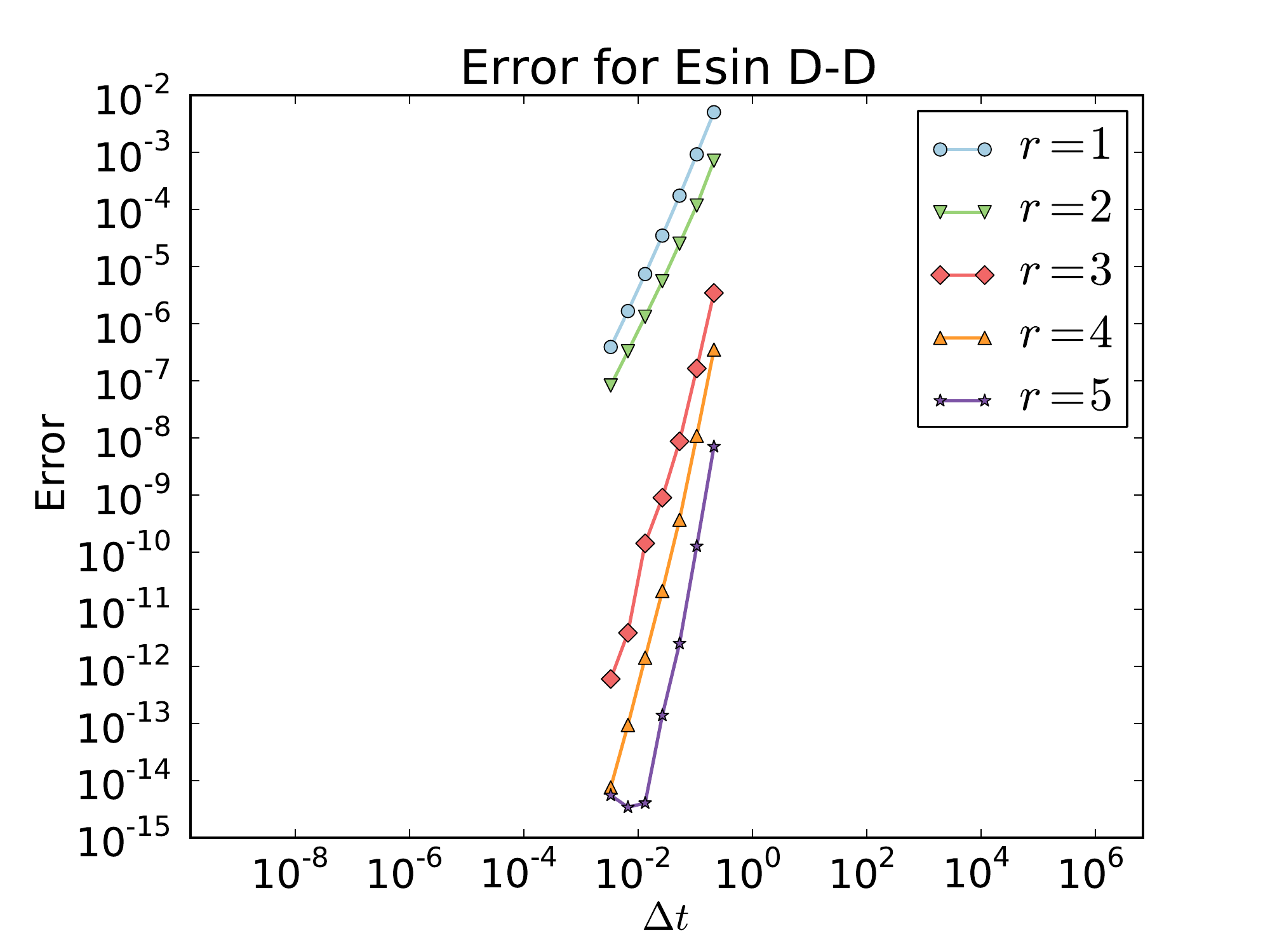}
   \includegraphics[width=\fw]{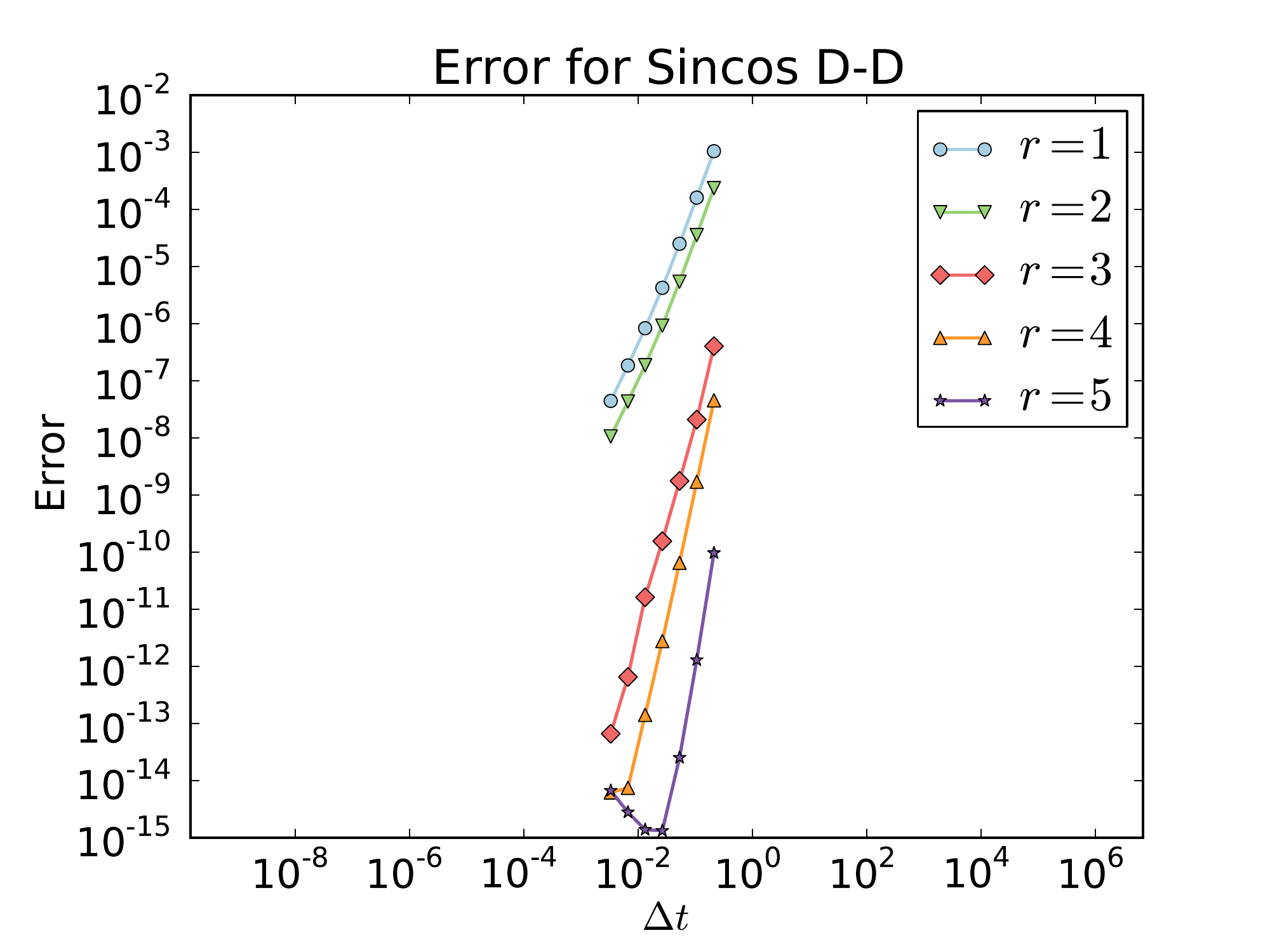}
  \includegraphics[width=\fw]{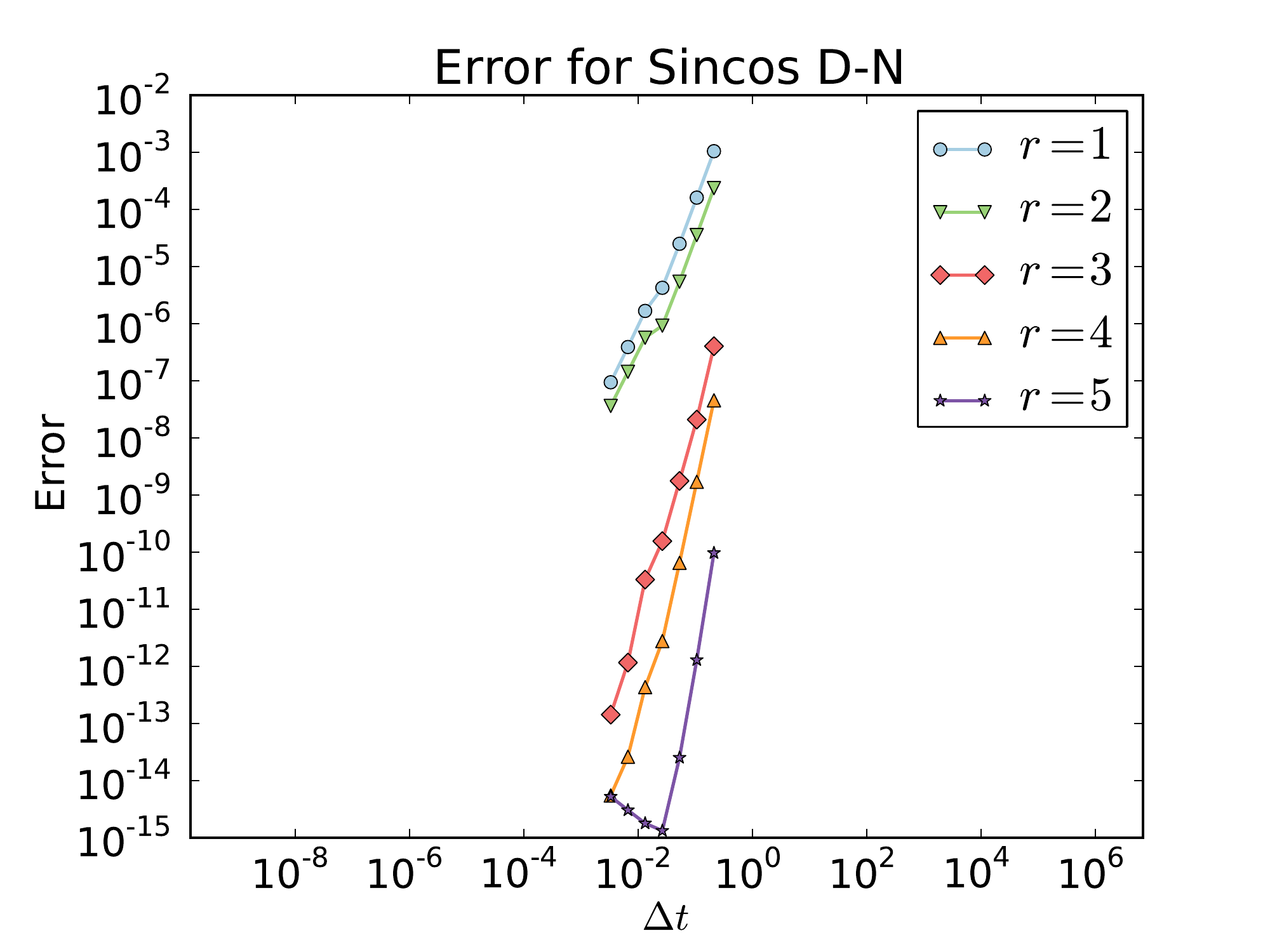}
   \includegraphics[width=\fw]{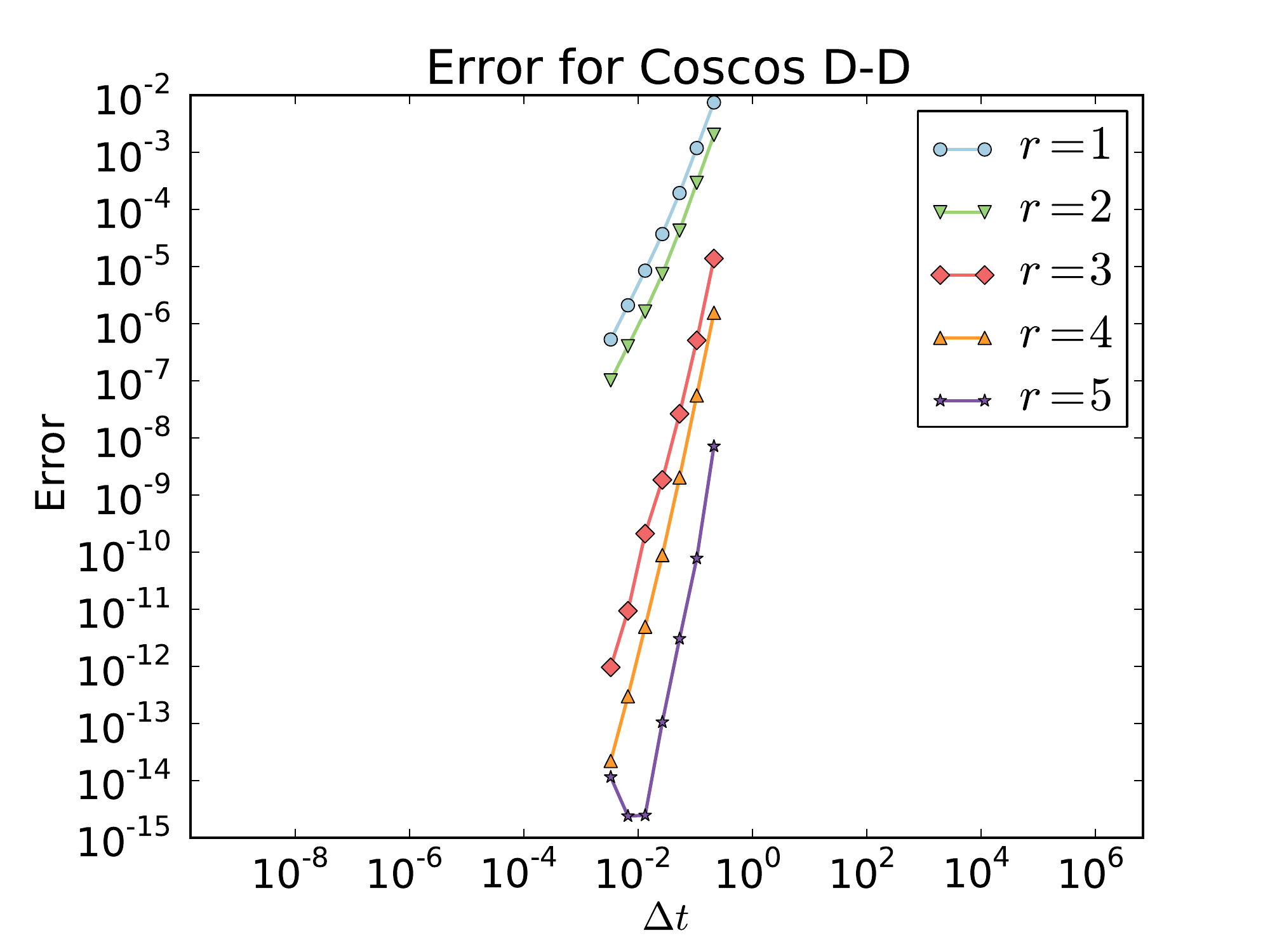}
   \includegraphics[width=\fw]{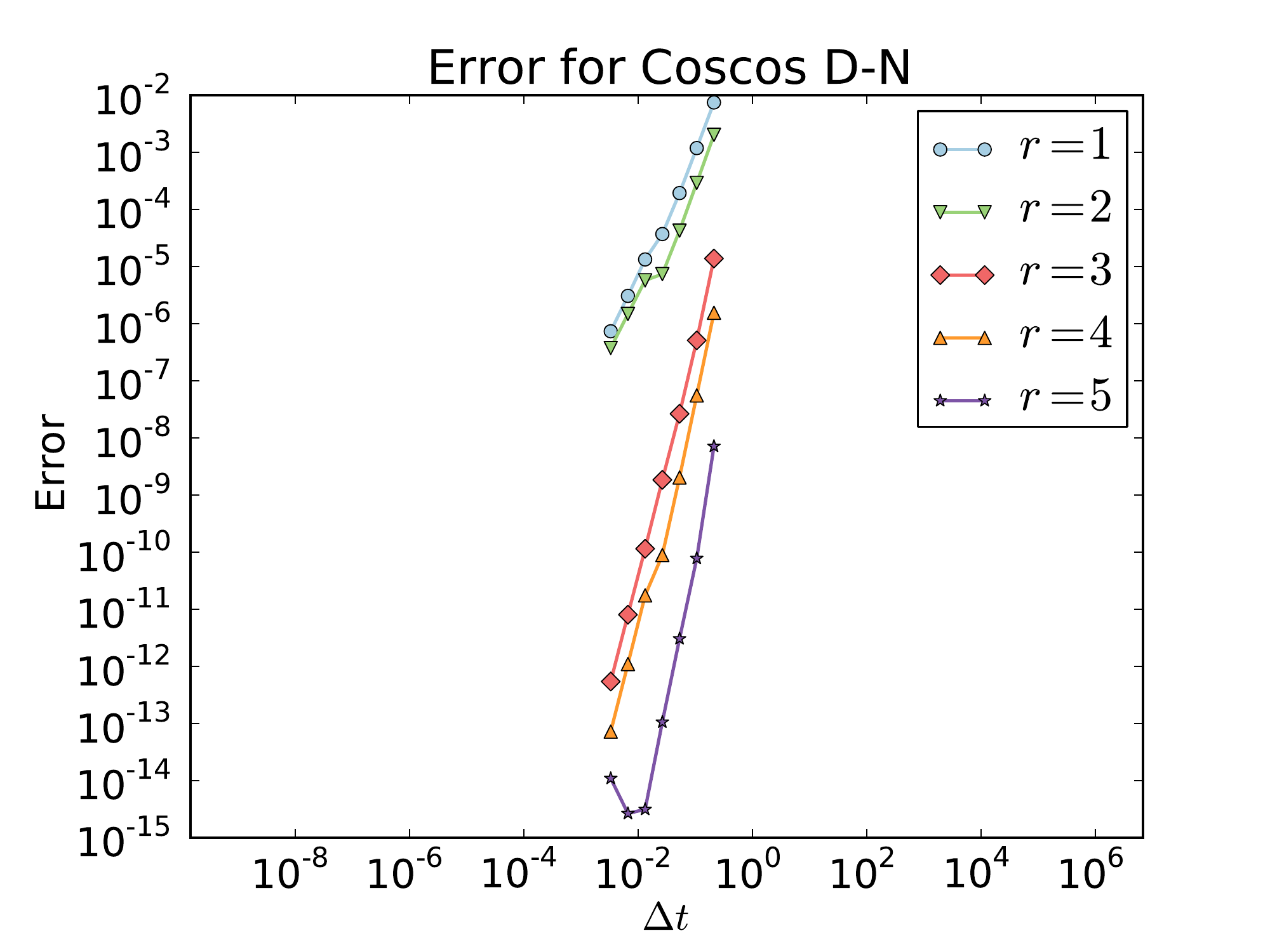}
   \includegraphics[width=\fw]{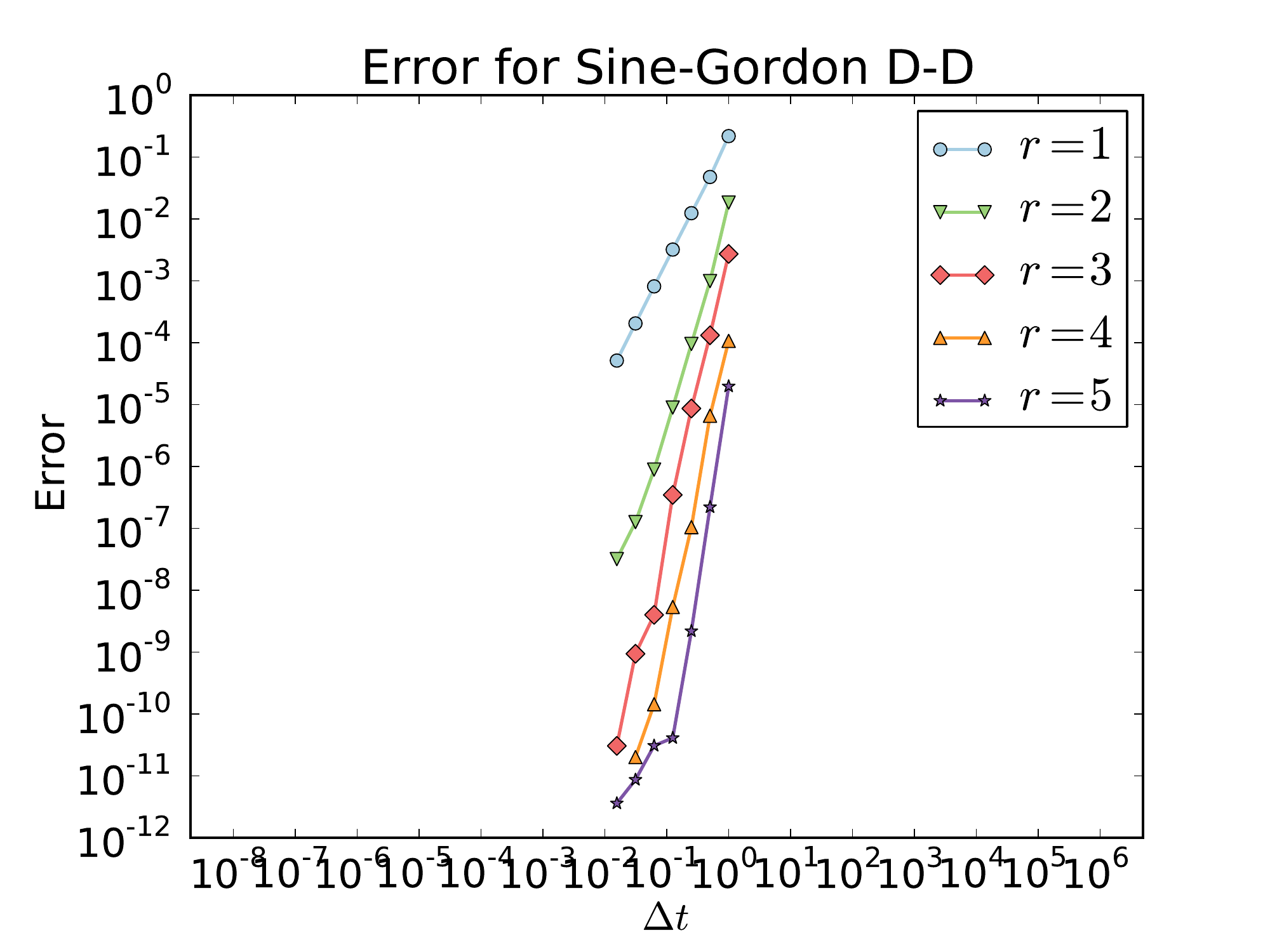}
   \caption{The error of the diamond scheme with varying $r$ applied
      to the Esin DD (Dirichlet/Dirichlet boundary conditions), Sincos DD,  Sincos DN (Dirichlet/Neumann), Coscos DD, Coscos DN and sine-Gordon DD  problems (see
      Table~\ref{table:samplednproblems}).  The Courant number
      is fixed at $\frac{1}{2}$ as $\Delta t$ is decreased.
      Table~\ref{table:erroroffourdnproblems} summarizes the observed convergence order
      given by the slope of these lines.  }
   \label{fig:error:init=diamond:physics=d_d_expasin:r=12345:s=2}
\end{figure}

Because
the exact solution is known for all the sample problems it is easy to
impose whatever boundary condition are desired on any spatial domain.  The
domains were chosen so the solutions were not periodic or symmetric in
any way, because while testing other potential methods it became
apparent that using `easy' problems gave false confidence in the
method.  To ensure only one thing was tested at a time, the
initialization scheme used was the diamond method.  As a comparison,
for one problem (Coscos DN), the boundary initialization was used.
Because the domains are smaller than the periodic
tests, a smaller number of diamonds was used, $N=2, 4, \ldots, 128$.
Figure~\ref{fig:error:init=diamond:physics=d_d_expasin:r=12345:s=2}
shows the error for the various test problems as $\Delta t$ is
decreased.  From this data the observed order of convergence is summarized in Table~\ref{table:erroroffourdnproblems}.

\begin{table}
\setlength{\tabcolsep}{1mm}
   \centering
   \begin{tabular}{c c c c c c c}
      \hline
      $r$ & \multicolumn{6}{c}{Order} \\
          & Esin DD & Sincos DD & Sincos DN & Coscos DD & Coscos DN & sine-Gordon DD \\ \hline
1 & 2.3 &  2.4 &  2.2 &  2.3 &  2.2 &  2.0 \\
2 & 2.2 &  2.4 &  2.1 &  2.4 &  2.1 &  3.2 \\
3 & 3.7 &  3.8 &  3.6 &  4.0 &  4.1 &  4.4 \\
4 & 4.2 &  4.5 &  3.8 &  4.3 &  4.1 &  4.5 \\
5 & 5.2 &  5.4 &  5.4 &  5.4 &  5.3 &  6.3 \\ \hline
   \end{tabular}
   \caption{Observed convergence orders for the initial--boundary value problems given in
      Table~\ref{table:samplednproblems}.  To one significant
      figure the order appears to be $r$, although in some cases it
   exceeds this.}      
   \label{table:erroroffourdnproblems}
\end{table}

\section{Conclusions} \label{chapter:conclusion} \label{sec:discussion}

The novel diamond mesh introduced in \cite{mclachlan2015multisymplectic} raised hopes that it would be suitable
for parallelisation and for a wide variety of initial--boundary value problems, with
local boundary closures not affecting the interior part of the scheme. We have
presented numerical evidence that this is indeed possible. It
still remains to prove stability and the observed orders of convergence.
That this may be possible is suggested by the known results
that the Reich method (for a wide class of equations of the form
\eqref{eq:hampde}, \eqref{eq:hampdemod}) has convergence order
at least $r$, and that it also preserves discrete forms
of arbitrary quadratic conservation laws, whenever
the PDE has any such \cite{sun2007}.

\begin{acknowledgements}
This research was supported by the Marsden Fund of the Royal Society Te Ap$\bar{\rm a}$rangi and by a Massey University PhD Scholarship.
\end{acknowledgements}

\bibliographystyle{siam}
\bibliography{refs}

\end{document}